\documentclass[11pt]{article}
\usepackage{amsmath,amsfonts,amssymb}
\usepackage{fancyhdr,graphicx,amsmath, amscd}
\usepackage{color}

\def\be{\begin{equation}}
\def\ee{\end{equation}}
\def\bea{\begin{eqnarray}}
\def\eea{\end{eqnarray}}
\def\bt{\begin{theorem}}
\def\et{\end{theorem}}
\def\bl{\begin{lemma}}
\def\el{\end{lemma}}
\def\br{\begin{remark}}
\def\er{\end{remark}}
\def\bc{\begin{corollary}}
\def\ec{\end{corollary}}
\def\bd{\begin{definition}}
\def\ed{\end{definition}}

\def\b{\beta}

\def\l{\lambda}
\def\m{\mu}

\def\r{\rho}
\def\s{\sigma}

\def\bbR{\mathbb{R}}

\def\b1{B_{1}^}

\def\ba{\begin{array}}
\def\ea{\end{array}}
\def\ben{\begin{enumerate}}
\def\een{\end{enumerate}}
\newtheorem{theorem}{Theorem}[section]
\newtheorem{lemma}{Lemma}[section]
\newtheorem{remark}{Remark}[section]
\newtheorem{proposition}{Proposition}[section]
\newtheorem{corollary}{Corollary}[section]

\newtheorem{definition}{Definition}[section]
\begin{document}
\thispagestyle{fancy}
\title{Inequalities for general mixed affine surface areas
}

\author{Deping Ye }
\date{}

\maketitle
\begin{abstract} Several general mixed affine surface areas are
introduced. We prove some important properties, such as, affine
invariance, for these general mixed affine surface areas. We also
establish new Alexandrov-Fenchel type inequalities,
Santal\'{o}-type inequalities, and affine isoperimetric
inequalities for these general mixed affine surface areas.

\end{abstract}

\section{Introduction}

There has been a growing body of work in isoperimetric
inequalities. The classical isoperimetric inequality, which
compares the surface area in terms of the volume, is an extremely
powerful tool in geometry and related areas. Relatively more
important results in the family of isoperimetric inequalities,
e.g., the celebrate Blaschke-Santal\'{o} inequality, have the
``affine invariant" flavor. These affine isoperimetric
inequalities compare two functionals associated with convex bodies
(or more general sets) where the ratio of the functionals is
invariant under non-degenerate linear transformations. Important
functionals include but are not limited to, the volume, $L_p$
affine surface areas, mixed $p$-affine surface areas, and general
affine surface areas.

\vskip 1mm The study of affine surface areas has a long history.
The notion of the classical affine surface area was first
introduced by Blaschke in 1923 \cite{Bl1}, and was first
generalized to the $L_p$ affine surface area for $p>1$ by Lutwak
in \cite{Lu1}. Since then, considerable attention has been paid to
the $L_p$ affine surface area, which is now at the core of the
rapidly developing $L_p$-Brunn-Minkowski theory \cite{CW1, CW2,
FGP, HMS, Lud2, LR1, Lut1993, LuYZh1, LYZ2004, SA1, SA2, TW1}
among others. The $L_p$ affine surface area was further extended
to all $p\in \bbR$ via geometric interpretations \cite{MW2, SW4,
SW5, Werner2007}. For a sufficiently smooth convex body $K$ in
$\bbR^n$, the $L_p$ affine surface area $as_{p}(K)$ of $K$ was
defined as in \cite{Lu1} ($p
>1$) and \cite{SW5} ($p <1$) by
\begin{equation*}
as_{p}(K)=\int
_{S^{n-1}}\bigg[h_{K}(u)^{1-p}f_{K}(u)\bigg]^{\frac{n}{n+p}}\,d\s
(u).
\end{equation*} Here $S^{n-1}$ is the boundary of the unit Euclidean ball $B^n_2$
in  $\mathbb{R}^n$, $\s$ is the usual surface area measure on
$S^{n-1}$, $h_{K}(u)$ is the support function of the convex body
$K$ at $u\in S^{n-1}$, and $f_{K}(u)$ is the curvature function of
$K$ at $u$, i.e., the reciprocal of the Gauss curvature $\kappa
_K(x)$ at the point $x\! \in \! \partial K$, the boundary of $K$,
that has $u$ as its outer normal. The $L_p$ affine surface area is
the key ingredient in many problems, such as, approximation of
convex bodies by polytopes \cite{Gr2, LuSchW, SW5},  theory of
valuation (see e.g. \cite{A1, A2, Klain1997, LudR}), and the $L_p$
affine isoperimetric inequality \cite{Lu1, WY}. Recently, Paouris
and Werner \cite{Paouris2010} linked the $L_p$ affine surface area
with the relative entropy of cone measures of $K$ and of its polar
$K^{\circ}=\{y\in \bbR^n: \langle x,y \rangle \leq 1, \forall x\in
K\}$, where $\langle x, y\rangle $ is the inner product of $x$ and
$y$.

\vskip 1mm In literature, two generalizations of the $L_p$ affine
surface area are important: the mixed $p$-affine surface area
\cite{Lut1987, Lu1, Wa2007, WernerYe2010} and the general affine
surface areas by Ludwig \cite{Ludwig2009, LR1}. The mixed
$p$-affine surface area, which involves $n$ convex bodies in
$\bbR^n$, takes the form
\begin{equation*}
as_p(K_1, \cdots, K_{n})=\int
_{S^{n-1}}\bigg[h_{K_1}(u)^{1-p}f_{K_1}(u)\cdots
h_{K_n}(u)^{1-p}f_{K_n}(u)\bigg]^{\frac{1}{n+p}}\,d\s (u).
\end{equation*} Clearly, $as_p(K)=as_p(K, \cdots, K)$ if all $K_i=K$. Moreover, the mixed
$p$-affine surface area contains many other important functionals
of convex bodies as special cases, such as, the usual surface area
and the dual mixed volume \cite{Lut1975}. The discovery of the
general affine surface areas owes to the valuation theory
\cite{LR1}. These general affine surface areas involve very
general functions, for instance, the $L_{\phi}$-affine surface
area associated with $\phi \in Conc(0,\infty)$ has the form
\begin{eqnarray*}as_{\phi} (K)&=&\!\!\!\!\int_{S^{n-1}}\!\!
\phi\bigg(\frac{1}{f_{K}(u)h_{K}^{n+1}(u)}\bigg)h_{K}(u)f_{K}(u)\,d\s
(u).
\end{eqnarray*} When $\phi(t)=t^{\frac{p}{n+p}}$ for $p>0$,
$as_{\phi}(K)=as_p(K)$. A fundamental result on the
$L_{\phi}$-affine surface area is the characterization theory of
upper-semicontinuous $SL(n)$ invariant valuation \cite{LR1}, that
is, every upper-semicontinuous, $SL(n)$ invariant valuation
vanishing on polytopes can be represented as a $L_{\phi}$-affine
surface area for some $\phi\in Conc(0,\infty)$. We refer readers
to \cite{Ludwig2009, LR1} for other general affine surfaces and
their properties.

\vskip 1mm In this paper, we introduce several general mixed
affine surface areas, which are the combinations of the above two
generalizations. Throughout the whole paper, $K \in C^2_+$ means
that $K$ has the origin in its interior, and has $C^2$ boundary
with everywhere strictly positive Gaussian curvature. Hereafter,
let $Conc(0, \infty)$ be the set of functions $\phi:\! (0,
\infty)\!\rightarrow\! (0,\infty)$ such that either $\phi$ is a
nonzero constant function, or $\phi$ is concave with $\lim
_{t\rightarrow 0}\phi(t)\!=\!0$ and $\lim_{t\rightarrow
\infty}\phi(t)/t\!=\!0$ (in this case, we set $\phi(0)\!=\!0$).
For all $\phi_i\!\in\! Conc(0,\infty)$ and all $K_i\!\in\! C^2_+$,
we define the general mixed $L_{\phi}$-affine surface area by
\begin{equation}\label{affine:conc}\!\! as(\phi_1,K_1;\!\cdots ;
\! \phi_n,K_n)\!=\!\!\int_{S^{n-1}}\!\! \prod
_{i=1}^n
\!\left[\phi_i\bigg(\!\frac{1}{f_{K_i}(u)h_{K_i}^{n+1}(u)}
\!\bigg)h_{K_i}(u)f_{K_i}(u)\right]^{\frac{1}{n}}\!\!\,d\s
(u).
\end{equation} Clearly, $as({\phi,\! K;\!
\cdots; \!\phi,\! K})\!=\!as_{\phi}(K)$, and $as(\phi,\! K_1;\!
\cdots;\! \phi,\! K_n)\!=\!as_p(K_1,\! \cdots,\! K_n)$ if
$\phi(t)\!=\!t^{\frac{p}{n+p}}$ for $p\!\geq\!0$. Hence, we
include the $L_p$ affine surface area for $p\!>\!0$ and the volume
(i.e., for $p=0$) as special cases. We show that the general mixed
$L_{\phi}$-affine surface area is affine invariant for all
$\phi_i\in Conc(0,\infty)$ (see Theorem \ref{affine:invariant:1}),
and also provide geometric interpretations of it (see Theorem
\ref{gemoertic:phi}). See Section 2 for other general mixed affine
surface areas and their properties, in particular, we prove a
duality result (Proposition \ref{dual:formula:star}).

\vskip 1mm The $L_p$ affine isoperimetric inequality for the $L_p$
affine surface area with all $p\in \bbR$  was proved in \cite{WY}
(see \cite{Lu1} for $p\geq 1$). The $L_p$ affine isoperimetric
inequalities are fundamental in many problems and have various
applications in, e.g., imaging recognition and PDE (e.g. \cite{Ga,
GaZ, Lu-O, Sch}). In particular, it was used by Andrews \cite{AN1,
AN2}, Sapiro and Tannenbaum \cite{ST1} to show the uniqueness of
self-similar solutions of the affine curvature flow and to study
its asymptotic behavior. Other related works include, e.g.,
\cite{CG, Haberl2009, LYZ1, LYZ2}. The $L_p$ affine isoperimetric
inequalities were extended to the mixed $p$-affine surface areas
\cite{Lu1, WernerYe2010} and the general affine surface areas
\cite{Ludwig2009}. In this paper, we prove analogous affine
isoperimetric inequalities for general mixed affine surface areas.
Let $B_{K_i}$ be the origin-symmetric (Euclidean) ball with
$|B_{K_i}|=|K_i|$ for all $i$.

\vskip 1mm \noindent{\bf Theorem
\ref{affine:isoperimetric:inequality}}  {\em Let all $K_i\in
C^2_+$ be convex bodies with centroid at the origin. \vskip 1mm
\noindent (i): If all $\phi_i\in Conc(0, \infty)$, then
$${as(\phi _1, K_1; \cdots ; \phi_n, K_n)}\leq {as(\phi_1, B_{K_1};
\cdots; \phi_n, B_{K_n})}.$$
Equality holds if all $K_i$ are ellipsoids that are dilates of one
another.

\vskip 1mm \noindent (ii) For all $\phi_i\in Conc(0,\infty)$ with
homogeneous degrees $r_i\in [0,1)$,
$$\left(\frac{as(\phi _1, K_1; \cdots ; \phi_n, K_n)}{as(\phi_1,
B^n_2; \cdots; \phi_n, B^n_2)}\right)^n \leq \prod _{i=1}^n
\left(\frac{|K_i|}{|B^n_2|}\right)^{1-2r_i},$$ with equality if
all $K_i$ are ellipsoids that are dilates of one another. }

\vskip 1mm Blaschke-Santal\'{o} inequality and the inverse
Santal\'{o} inequality are model examples of affine isoperimetric
inequalities. Blaschke-Santal\'{o} inequality was proved by the
$L_1$ affine isoperimetric inequality in \cite{Santalo1949} (where
the Santal\'{o} point instead of the centroid was used). Note that
if $K$ has Santal\'{o} point at the origin, then its polar body
$K^\circ$ has centroid at the origin (see \cite{Sch} for details).
The inverse Santal\'{o} inequality was due to Bourgain and Milman
\cite{BM} (see also \cite{GK2, MilmanPajor2000, Nazarov2009,
Pisier1989}). They were successfully extended to affine surface
areas in \cite{Lu1, Wa2007, WY, WernerYe2010} among others. Here
we prove the following Santal\'{o}-type inequalities for the
general mixed $L_{\phi}$-affine surface areas.

\vskip 1mm \noindent {\bf Theorem \ref{Stantalo:phi}} {\em Let
$K_i\in C^2_+$ and $\phi_i \in Conc(0,\infty)$ be homogeneous of
degrees $r_i\in [0,1)$. Then
$$as(\phi_1, K_1; \cdots; \phi_n, K_n)as(\phi_1, K_1^\circ; \cdots
; \phi_n, K_n^\circ)\leq n^2 \prod_{i=1}^n \big[\phi_i(1)^2\
|K_i|\ |K_i^\circ|\big]^{\frac{1}{n}}.$$ Moreover, if all $K_i\in
C^2_+$ have centroid at the origin, then
$$as(\phi_1, K_1; \cdots; \phi_n, K_n)as(\phi_1, K_1^\circ;
\cdots ; \phi_n, K_n^\circ)\leq \big[as(\phi_1, B^n_2; \cdots;
\phi_n, B^n_2)\big]^2,$$ with equality if all $K_i$ are
origin-symmetric ellipsoids that are dilates of one another.}

\vskip 1mm Another powerful inequality in geometry is the
well-known classical Alexandrov-Fenchel inequality for mixed
volume (see \cite{ Ale1937, Bus1958, Sch}). Here we show
Alexandrov-Fenchel type inequalities for the general mixed affine
surface areas, which have similar forms to the classical
Alexandrov-Fenchel inequality. However, they do not imply the
classical Alexandrov-Fenchel inequality. See also \cite{Lut1975,
Lut1987, WernerYe2010} for more such type inequalities.

\vskip 1mm \noindent {\bf Proposition \ref{inequality:mixed:conc}}
{\em Let all $K_i\in C^2_+$ and $\phi_i\in Conc(0,\infty)$. Then,
for $1 \leq m \leq n$,
$$ as^m(\phi_1,\! K_1;\! \cdots\! ;\! \phi_n,\! K_{n})\! \leq \!\!\prod_{i=0}^
{m-1}as(\phi_1,\! K_1;\! \cdots;\! \phi_{n-m},\! K_{n-m};
\underbrace{\phi_{n-i},\! K_{n-i};\! \cdots;\! \phi_{n-i},\!
K_{n-i}}_m).
$$ Equality holds if (1) all $K_i$ coincide and $\phi_i\!=\!\l_i\phi_n$
for some $\l_i\!>\!0$, $i\!=\!n-m+1,\! \cdots\!, n$, or (2)
$\phi_i=\l_i\phi_n$, $K_i=\eta _i K_n$ for some $\l_i, \eta_i>0$,
$i=n-m+1, \cdots, n$, and $\phi_n$ is homogeneous of degree $r\in
[0,1)$. The equality holds trivially if $m=1$.

\vskip 1mm \noindent In particular, $as^n(\phi_1,K_1; \cdots;
\phi_n, K_{n}) \leq  as_{\phi_1}(K_1) \cdots as_{\phi_n}(K_n). $ }

\vskip 1mm This paper is organized as follows. In Section 2, we
introduce several general mixed affine surface areas. We provide
geometric interpretations of them and prove some important
properties of them, such as, affine invariant properties. In
Section 3, we establish new Alexandrov-Fenchel type inequalities,
Santal\'{o}-type inequalities, and affine isoperimetric
inequalities for these general mixed affine surface areas. Section
4 dedicates to the general $i$-th mixed affine surface areas.
Similar Santal\'{o}-type and affine isoperimetric inequalities are
also proved.


\section{General mixed affine surface areas}


\subsection{General mixed $L_{\phi}$- and $L_{\psi}$-affine surface areas}

\vskip 1mm The $L_p$ affine surface area for $-n<p<0$ is
associated with the $L_{\psi}$-affine surface area, where $\psi\in
Conv(0,\infty)$ in \cite{Ludwig2009}. Hereafter, $Conv(0,\infty)$
is the set of functions $\psi: (0,\infty)\rightarrow (0,\infty)$
such that either $\psi$ is a nonzero constant function, or $\psi$
is convex with $\lim_{t\rightarrow 0} \psi(t)=\infty$ and $
\lim_{t\rightarrow \infty} \psi(t)=0$ (in this case, we set
$\psi(0)=\infty$). For $\psi _i \in Conv(0,\infty)$ and $K_i\in
C^2_+$, we define the general mixed $L_{\psi}$-affine surface area
as
\begin{equation*} \!\!as(\psi_1,\!K_1;\!\cdots;\! \psi_n,\!K_n)\!=\!\int_{S^{n-1}}\!\!\prod
_{i=1}^n\!\!
\left[\!\psi_i\bigg(\frac{1}{f_{K_i}(u)h_{K_i}^{n+1}(u)}\bigg)h_{K_i}(u)f_{K_i}(u)\!\right]^{\frac{1}{n}}\!\!\,d\s
(u). \end{equation*}  In particular, $as(\psi,K;\cdots;
\psi,K)\!=\!as_{\psi}(K)$ is the $L_{\psi}$-affine surface area of
$K$ introduced in \cite{Ludwig2009}.  When all
$\psi_i(t)\!=\!t^{\frac{p}{n+p}}$ for $-n\!<\!p\!\leq\!0$, one
gets the mixed $p$-affine surface area for $-n\!<\!p\!\leq\!0$. In
particular, one includes the $L_p$ affine surface area for
$-n\!<\!p\!<\!0$ as a special case. If all $\psi_i\!=\!\psi$, then
we use $as_{\psi}(K_1\!, \cdots\!, K_n)$ to represent
$as(\psi\!,K_1;\! \cdots\!; \psi,\! K_n)$. We use $as(\psi_1,\!
\cdots\!, \psi_n; K)$ instead of $as(\psi_1,\! K;\! \cdots\!;
\psi_n, \!K )$ when all $K_i\!=\!K$. Clearly, $as(\psi_1,\!
\cdots\!, \psi_n;\! B^n_2)\!=\!\big[\psi_1(1)\!\cdots
\!\psi_n(1)\big]^{\frac{1}{n}} n |B^n_2|.$

\vskip 1mm We always assume that $\phi\!\! \in\!\! Conc(0,\infty)$
is nonzero. As above, if all $\phi_i\!\!=\!\!\phi$, we write
$as_{\phi}(K_1,\! \cdots\!, K_n)$ for $as(\phi,\!K_1;\! \cdots\!;
\phi,\! K_n)$. We use $as(\phi_1,\! \cdots\!, \phi_n;\! K)$ for
$as(\phi_1,\! K;\! \cdots\!; \phi_n,\! K)$ if all $K_i\!=\!K$.
Clearly, $as(\phi_1,\! \cdots\!, \phi_n;\!
B^n_2)\!=\!\big[\phi_1(1)\!\cdots\! \phi_n(1)\big]^{\frac{1}{n}} n
|B^n_2|.$


\vskip 1mm The following theorem gives a geometric interpretation
for the general mixed $L_{\phi}$-affine surface area by the
illumination surface body. Similar geometric interpretations can
also be obtained by the surface body \cite{SW5, WY}.

\bd {\rm ({\bf Illumination surface body}) \cite{WernerYe2010}}
Let $s\!\geq\! 0$ and $f:\!\partial K\!\rightarrow\! \bbR$ be a
nonnegative, integrable function. The illumination surface body
$K^{f,s}$ is defined as
$$\displaystyle K^{f,s}=\left\{x: \int _{\partial K \cap
\overline{[x,K]\backslash K}}f\,d\mu_K\leq s\right\}.$$ Here,
$\mu_K$ is the usual surface measure of $\partial K$, $[x,K]$
denotes the convex hull of $x$ and $K$, $A\setminus B=\{z: z\in A,
\ but \ z\notin B\}$, and $\bar{A}$ is the closure of $A$. \ed

 \bt\label{gemoertic:phi} Let $K, K_i\in C^2_+$ and $\phi_i\in Conc(0,\infty)$, $i=1, \cdots,
 n$. Let $f: \partial K\rightarrow \bbR$ be the function
\begin{equation}\label{function}
f(N_K^{-1}(u))=f_K(u)^{\frac{n-2}{2}}\prod_{i=1}^n\left[\phi_i\bigg(\!\frac{1}{f_{K_i}(u)h_{K_i}^{n+1}(u)}
\!\bigg)h_{K_i}(u)f_{K_i}(u)\right]^{\frac{1-n}{2n}},
\end{equation} where $N_K^{-1}: S^{n-1}\rightarrow \partial K$ is
the inverse Gauss map, that is, $N^{-1}_K(N_K(x))=x$ for all $x\in
\partial K$. Let $c_n=2|B^{n-1}_2|^{\frac{2}{n-1}}.$ Then,
$$as(\phi_1, K_1; \cdots ; \phi_n, K_n)=\lim_{s\rightarrow 0} c_n\
\frac{|K^{f,s}|-|K|}{s^{\frac{2}{n-1}}}.$$
 \et
\noindent {\bf Proof.} Theorem 4.1 and its following remark in
\cite{WernerYe2010} imply that \begin{eqnarray*} \lim
_{s\rightarrow 0} c_n \frac{|K^{f,s}|-|K|}{
s^{\frac{2}{n-1}}}&=&\int_{S^{n-1}}
\frac{f_K(u)^{\frac{n-2}{n-1}}}{f(N_K^{-1}(u))^{\frac{2}{n-1}}}\,d\s
(u)\\ & =&\int _{S^{n-1}} \prod _{i=1}^n
\left[\phi_i\bigg(\!\frac{1}{f_{K_i}(u)h_{K_i}^{n+1}(u)}\!
\bigg)h_{K_i}(u)f_{K_i}(u)\right]^{\frac{1}{n}}\,d\s(u) \\
&=& as(\phi_1, K_1; \cdots ; \phi_n, K_n),
\end{eqnarray*} where the second equality is by (\ref{function}) and the last equality is by
(\ref{affine:conc}).

\vskip 1mm \noindent {\bf Remark.} Replacing $\phi_i\in
Conc(0,\infty)$ by $\psi _i\in Conv(0,\infty)$, one gets the
geometric interpretation of the general mixed $L_{\psi}$-affine
surface area. That is,
$$as(\psi_1, K_1; \cdots ; \psi_n, K_n)=\lim_{s\rightarrow 0} c_n\
\frac{|K^{\tilde{f},s}|-|K|}{s^{\frac{2}{n-1}}},$$ where the
function $\tilde{f}$ takes the form
$$\tilde{f}(N_K^{-1}(u))=f_K(u)^{\frac{n-2}{2}}\prod_{i=1}^n
\left[\psi_i\bigg(\!\frac{1}{f_{K_i}(u)h_{K_i}^{n+1}(u)}\!\bigg)h_{K_i}(u)f_{K_i}(u)\right]^{\frac{1-n}{2n}}.$$

\vskip 1mm A function $\phi\in Conc(0,\infty)$ is homogeneous of
degree $r$ if $\phi(\l t)=\l ^r\phi(t)$ for all $\l>0$ and $t>0$.
This further implies that $\phi(t)=\phi(1) t^r$ with $r\in [0,1)$
and $\phi(t)\phi(1/t)=\phi(1)^2$. Similarly, a function $\psi \in
Conv(0,\infty)$ is homogeneous of degree $r$ if $\psi (\l t)=\l ^r
\psi(t)$ for all $\l, t>0$. This further implies that
$\psi(t)=\psi (1) t^r$ with $r\in (-\infty, 0]$ and
$\psi(t)\psi(1/t)=\psi(1)^2$.

\vskip 1mm The following theorem gives the affine invariant
property for the general mixed $L_{\phi}$-affine surface area.
This result was proved in \cite{Ludwig2009} for all $\phi_i=\phi$
and $K_i=K$.
 \bt
\label{affine:invariant:1}  Let $T: \bbR^n \rightarrow\bbR^n $ be
an invertible linear transform. For all $\phi_i\in Conc(0,\infty)$
and $K_i\in C^2_+$,
$$as(\phi_1, TK_1; \cdots ; \phi_n, TK_{n})=as(\phi_1, K_1; \cdots ; \phi_n, K_{n}), \ for \ \ \  |det(T)|=1.$$

\noindent If in addition, $\phi_i\in Conc(0,\infty)$ are
homogeneous of degrees $r_i\in [0,1)$ for $i=1, \cdots, n$, that
is, $\phi_i(t)=\phi_i(1) t^{r_i}$, then
$$as(\phi_1, TK_1; \cdots ; \phi_n, TK_{n})=|det(T)|^{1-\frac{2\sum _{i=1}^n r_i}{n}} as(\phi_1,
K_1; \cdots ; \phi_n, K_{n}).$$ \et

\vskip 1mm \noindent {\bf Remark.} Replacing $\phi_i\in
Conc(0,\infty)$ by $\psi_i\in Conv(0,\infty)$, one has the affine
invariant property for the general mixed $L_{\psi}$-affine surface
area; namely, $$as(\psi_1, TK_1; \cdots ; \psi_n,
TK_{n})=as(\psi_1, K_1; \cdots ; \psi_n, K_{n}),\ \ for \ \ \
|det(T)|=1.$$ If in addition, $\psi_i(t)=\psi_i (1) t^{r_i}$, then
$$as(\psi_1, TK_1; \cdots ; \psi_n, TK_{n})=|det(T)|^{1-\frac{2\sum
_{i=1}^n r_i}{n}} as(\psi_1, K_1; \cdots ; \psi_n, K_{n}).$$

\vskip 1mm \noindent {\bf Proof.} Lemma 12 of \cite{SW5} implies
that, for all $u\in S^{n-1}$,
\begin{eqnarray}\label{Curvature:K:TK}
f_K(u)=\frac{f_{TK}\left(v\right)}{\  det ^2(T) \
\|T^{-1t}(u)\|^{n+1} },
\end{eqnarray} where $v=\frac{T^{-1t}(u)}{\|T^{-1t}(u)\|}\in S^{n-1}$ with
$\|\cdot\|$ standing for the Euclidean norm, $A^t$ denotes the
usual adjoint of $A$, and $A^{-1}$ is the inverse of $A$ for an
operator $A$. On the other hand,
\begin{equation}\label{support:K:TK}h_K(u)=\|T^{-1t}(u)\|\
h_{TK}(v).\end{equation} Thus, we have
\begin{eqnarray} && {f_{TK}(v)h_{TK}^{n+1}(v)}=det^2(T)\
f_K(u)h_K^{n+1}(u), \label{Curvature:K:TK:2}
\\&& {f_{TK}(v)h_{TK}(v)}=det^2(T)\|T^{-1t}(u)\|^n
f_K(u)h_K(u).\label{Curvature:K:TK:3}
\end{eqnarray} Lemma 10 in
\cite{SW5} implies that, for any integrable function $g:
S^{n-1}\rightarrow \bbR$,
\begin{eqnarray*} \int_{S^{n-1}}g(v)
\,d\s(v)=\frac{1}{|{det}(T)|}\int_{S^{n-1}}g\left(\frac{T^{-1t}(u)}{\|T^{-1t}(u)\|}\right)
\|T^{-1t}(u)\|^{-n}\,d\s(u).
\end{eqnarray*} Hence, by (\ref{Curvature:K:TK:2}) and
(\ref{Curvature:K:TK:3}), one has
\begin{eqnarray*}
as(\phi_1,\! TK_1;\! \cdots\! ; \phi_n,\!
TK_n)\!\!=\!\!\!\int_{S^{n-1}}\!\! \prod _{i=1}^n
\!\left[\phi_i\bigg(\!\frac{1}{f_{TK_i}(v)h_{TK_i}^{n+1}(v)}\!\bigg)h_{TK_i}(v)f_{TK_i}(v)\right]^{\frac{1}{n}}\!\!\,d\s
(v)\\
=|det(T)|\int_{S^{n-1}}\!\! \prod _{i=1}^n
\!\left[\phi_i\bigg(\!\frac{1}{det(T)^2
f_{K_i}(u)h_{K_i}^{n+1}(u)}\!\bigg)h_{K_i}(u)f_{K_i}(u)\right]^{\frac{1}{n}}\!\!\,d\s
(u).
\end{eqnarray*}

\vskip 1mm \noindent Clearly, $|det(T)|=1$ implies $as(\phi_1,
TK_1; \cdots ; \phi_n, TK_n)=as(\phi_1, K_1; \cdots ; \phi_n,
K_n).$

\vskip 1mm  If in addition, $\phi_i$ are homogeneous of degrees
$r_i\in [0,1)$ for $i=1, \cdots, n$, then
$$\phi_i\bigg(\!\frac{1}{det(T)^2
f_{K_i}(u)h_{K_i}^{n+1}(u)}\!\bigg)=det(T)^{-2r_i} \phi
_i\bigg(\!\frac{1}{f_{K_i}(u)h_{K_i}^{n+1}(u)}\!\bigg),\ \forall
u\in S^{n-1},
$$ and therefore,
\begin{eqnarray*}
as(\phi_1, TK_1; \cdots ; \phi_n, TK_n)= |det(T)|^{1-\frac{2\sum
_{i=1}^n r_i}{n}} as(\phi_1, K_1; \cdots; \phi_n, K_n).
\end{eqnarray*}


\subsection{General mixed $L_{\phi}^*$- and $L_{\psi}^*$-affine surface areas}
\vskip 1mm In \cite{Ludwig2009}, Ludwig showed that the $L_p$
affine surface area for $p\!<\!-n$ is associated with the
$L_{\psi}^*$-affine surface area for $\psi\!\in\! Conv(0,\infty)$.
Here, we define the general mixed $L_{\psi}^*$-affine surface area
by
\begin{equation*} as^*(\psi_1, K_1; \cdots ; \psi_n, K_n)=\int _{S^{n-1}} \prod
_{i=1}^n
\left[\frac{\psi_i({f_{K_i}(u)h_{K_i}^{n+1}(u)})}{h_{K_i}(u)^n}\right]^{\frac{1}{n}}\,d\s
(u),
\end{equation*} for $\psi_i\!\in\! Conv(0,\infty)$
and $K_i\!\in\! C^2_+$, $i\!=\!1,\cdots, n$. As above,
$as^*_{\psi}(K_1,\!\cdots\!, K_n)\!=\!as^*(\psi,\!K_1;\!\cdots\!;
\psi,\! K_n)$ and $as_{\psi}^*(K)\!=\!as_{\psi}^*(K,\! \cdots,\!
K)$. If all $K_i$ coincide with $K$, we use $as^*(\psi_1,\!
\cdots\!, \psi _n;\! K)$ instead of $as^*(\psi_1,\! K;\! \cdots\!
; \psi_n,\! K)$. If $\psi(t)\!=\!t^{\frac{n}{n+p}}$ for
$p\!<\!-n$, then $as_{\psi}^*(K_1,\! \cdots,\! K_n)\!=\!as_p(K_1,
\!\cdots,\! K_n)$ and hence, $as_{\psi}^* (K)$ is the $L_p$ affine
surface area for $p\!<\!-n$. In particular, $as^*(\psi_1, \cdots,
\psi_n; B^n_2)\!=\!\big[\psi_1(1)\!\cdots\!
\psi_n(1)\big]^{\frac{1}{n}} n |B^n_2|.$

\vskip 1mm The following proposition gives the duality relation
between $L_{\psi}$- and $L_{\psi}^*$-affine surface areas.

\vskip 1mm  \begin{proposition}\label{Duality:psi} Let all $K_i\in
C^2_+$ be convex bodies, such that, $K_i=\l_{i}K$ for some convex
body $K\in C^2_+$ and $\l_i>0$, $i=1, \cdots, n$. For all
$\psi_i\in Conv(0,\infty)$,
$$as^*(\psi_1, K_1; \cdots ; \psi_n, K_n)=as(\psi_1,
K_1^\circ;\cdots ; \psi_n, K_n^\circ ).$$ In particular,
$as^*(\psi_1, \cdots, \psi_n; K)=as(\psi_1, \cdots, \psi_n;
K^\circ).$ \end{proposition}

\vskip 1mm \noindent {\bf Proof.} Define $y:\!
S^{n-1}\!\rightarrow \!
\partial K^{\circ}$ by $y(u)\!=\!\rho_{K^\circ}(u)u$ with
$\rho_{K^\circ}(u)$ the radius function of $K^\circ $ at the
direction $u$, that is, $\rho_{K^\circ}(u)=\max\{\l>0: \l u\in
K^\circ \}.$ Note that $h_K(u)\rho_{K^\circ}(u)=1$ for all
directions $u\in S^{n-1}$. The Jacobian $Jy$ is (see e.g.
\cite{Hug} ) $$Jy(u)=\frac{\r _{K^{\circ}}(u)^{n-1} }{\langle u,
N_{K^{\circ}}(\r _{K^{\circ}}(u)u)\rangle}, \ \ \ a.s. \ \ on \ \
S^{n-1}. $$ The area formula (see e.g.\ \cite{Federer1969})
implies that for every a.s.\ defined function $g:\!
S^{n-1}\!\rightarrow\! [0,\infty]$, one has $\int _{S^{n-1}}
g(u)Jy(u)\,d\s(u)\!=\!\int _{\partial K^{\circ}}
g\big(\frac{y}{\|y\|}\big)\,d\mu _{K^{\circ}(y)}.$ Setting
$$g(u)=\prod _{i=1}^n \!\left[\psi_i\bigg( \frac{1}{\l _i^{2n}
f_{K}(u)h_{K}^{n+1}(u)}\!\!\bigg)h_{K}(u)f_{K}(u)\!\right]^{\frac{1}{n}}\frac{\langle
u,N_{K^\circ}(y(u)) \rangle }{\r _{K^\circ}(u)^{n-1}},$$ one has
\begin{eqnarray}\!\!\!\! &&as(\psi _1, K_1; \cdots ;
\psi _n, K_n)\! =\! as(\psi _1, \l_1 K; \cdots ; \psi _n, \l _n K )\nonumber\\
\!\!\!\! &&\!\!=\!\!\l_1 \cdots \l_n \ \int_{S^{n-1}} \prod
_{i=1}^n \!\left[\psi_i\bigg( \frac{1}{\l _i^{2n}
f_{K}(u)h_{K}^{n+1}(u)}\!\!\bigg)h_{K}(u)f_{K}(u)\!\right]^{\frac{1}{n}}\!\!\,d\s
(u)\nonumber\\
\!\!\! \!&&\!\!=\!\!\l_1\! \!\cdots\!\! \l_n \!\!
\int_{K^{\circ}}\!\! \prod _{i=1}^n \!\!\left[\!\psi_i\!\bigg(\!\!
\frac{\l _i^{-2n}}{
f_{K}\!\!\left(\!\!\frac{y}{\|y\|}\!\right)\!\!
h_{K}^{n+1}\!\!\left(\!\frac{y}{\|y\|}\!\!\right)}\!\bigg)h_{K}\!\!\left(\!\!\frac{y}{\|y\|}\!\!\right)\!\!
f_{K}\!\!\left(\!\!\frac{y}{\|y\|}\!\!\right)\!\!\right]^{\frac{1}{n}}\!\!\frac{\langle
{y}, N_{K^{\circ}}(y)\rangle}{\r _{K^{\circ}}(\frac{y}{\|y\|})^{n}
} \,d\m_{K^\circ} (y)\nonumber\\
\!\!\!&&\!\!=\!\!\l_1 \!\!\cdots\!\! \l_n \!\!
\int_{K^{\circ}}\!\! \prod _{i=1}^n \!\!\left[\!\psi_i\!\bigg(\!\!
\frac{\l _i^{-2n}}{
f_{K}\!\!\left(\!\!\frac{y}{\|y\|}\!\!\right)\!
h_{K}^{n+1}\!\left(\!\frac{y}{\|y\|}\!\right)\!}\!\bigg)\right]^{\frac{1}{n}}\!\!\!\!\!
h^{n+1}_{K}\!\!\left(\!\!\frac{y}{\|y\|}\!\!\right)\!\!
f_{K}\!\!\left(\!\!\frac{y}{\|y\|}\!\!\right)\!\! {\langle {y},
N_{K^{\circ}}(y)\rangle}\! \,d\m_{K^\circ} (y),\nonumber\\
\label{duality:integral}
\end{eqnarray} where the third equality is by $\|y\|\!=\!\r
_{K^{\circ}}(\frac{y}{\|y\|})$ and the last equality is by
$h_K(u)\!=\!\r _{K^{\circ}}(u)$ for all directions $u\in S^{n-1}$.

\vskip 1mm Now we let $v=N_{K^{\circ}}(y)$. Hug \cite{Hug} proved
that for almost all $y\in \partial K^{\circ}$,
$$h^{n+1}_{K} \left( \frac{y}{\|y\|} \right)
f_{K} \left( \frac{y}{\|y\|}
\right)=\frac{1}{h^{n+1}_{K^\circ}(v)f_{K^{\circ}}(v)}.$$
Combining with $ \,d\m_{K^\circ} (y)=f_{K^{\circ}}(v)\,d\s(v)$ and
by $(\l_i K)^\circ =\l _i^{-1} K^\circ$, the equality
(\ref{duality:integral}) equals to
$$\l_1 \cdots \l_n \int_{S^{n-1}} \prod _{i=1}^n  \left[ \psi_i
\bigg({\l
_i^{-2n}}{h^{n+1}_{K^\circ}(v)f_{K^{\circ}}(v)}\!\bigg)\right]^{\frac{1}{n}}
\ \frac{\,d\s(v)}{h^n_{K^{\circ}}(v)} =as ^*(\psi _1, K_1^\circ;
\cdots ; \psi_n, K_n^\circ ), $$  which completes the proof.

\vskip 1mm \noindent {\bf Remark.} When all $\psi_i\!=\!\psi$,
this result was proved in \cite{Ludwig2009}, i.e.,
$as_{\psi}^*(K)\!=\!as_{\psi}(K^\circ)$ for all $K\!\in\! C^2_+$
and $\psi\!\in\! Conv(0,\infty).$ In particular, if
$\psi(t)=t^{\frac{n}{n+p}}$ for $p<-n$, then
$as_p(K)\!=\!as_{\frac{n^2}{p}}(K^\circ)$ \cite{WY}. In general,
one {\em cannot} expect $as^*(\psi_1, K_1; \cdots ; \psi_n,
K_n)=as(\psi_1, K_1^\circ;\cdots ; \psi_n, K_n^\circ )$ even if
all $\psi_i$ are equal to some $\psi\!\in\! Conv(0, \infty)$ and
all $K_i$ are ellipsoids. To this end, let $n=2$ and
$\psi(t)=\frac{1}{t}$. For any $2$-dimensional convex body $K\in
C^2_+$,
\begin{eqnarray*}as_{\psi}(K, B^2_2)\!\! =\!\! \int_{S^{1}}\! \left[\psi
\bigg(\frac{1}{f_{K}(u)h_{K}^{3}(u)}\bigg)h_{K}(u)f_{K}(u)\right]^{\frac{1}{2}}\!\!\,d\s_1
(u)\!\! =\!\! \int_{S^{1}}\! h^2_{K}(u)f_{K}(u)\,d\s_1(u),
\end{eqnarray*} where $\s_1$ refers to the spherical measure of $S^1$.
On the other hand,
$$as^*_{\psi}(K^\circ, B^2_2)=\int _{S^{1}}
\frac{\,d\s_1(u)}{\sqrt{f_{K^\circ}(u)h^5_{K^\circ}(u)}}.$$ Now
let the (invertible) affine map $T: \bbR^2\rightarrow \bbR^2$ be
$T(x_1, x_2)=(x_1, 2x_2)$. Then, by formulas
(\ref{Curvature:K:TK}) and (\ref{support:K:TK}),
$$h_{TB^2_2}(u)=\|Tu\|, \ \ \ f_{TB^2_2}(u)=\frac{4}{\|Tu\|^3}, \ \ \forall u\in
S^1.$$ Therefore, one has
$$as_{\psi}(TB^2_2, B^2_2)\!\!=\!\!\int _{S^1}\!
f_{TB^2_2}(u)h^2_{TB^2_2}(u)\,d\s_1(u)\!\!=\!\!
\int_{S^1}\!\frac{4}{\|Tu\|}\,d\s_1(u)\!\!=\!\!\int
_{S^1}\frac{4}{\sqrt{1+3u_2^2}}\,d\s_1(u).$$ As
$(TB^2_2)^\circ=T^{-1t}B^2_2$, one has
$$h_{(TB^2_2)^{\circ}}(u)=\|T^{-1}u\|, \ \ \
f_{(TB^2_2)^\circ}(u)=\frac{4}{\|T^{-1}u\|^3}, \ \ \forall u\in
S^1,$$ and hence
$$as^*_{\psi}((TB^2_2)^\circ, B^2_2)=\int_{S^1}\frac{\,d\s_1(u)}
{\sqrt{4\|T^{-1}u\|^2}}=\int
_{S^1}\frac{1}{\sqrt{1+3u_1^2}}\,d\s_1(u).$$ Clearly
$as_{\psi}(TB^2_2, B^2_2)=4as^*_{\psi}((TB^2_2)^\circ, B^2_2)$ by
the rotational invariance of the spherical measure $\s_1$.

 \vskip 1mm For all $\phi_i\in Conc(0,\infty)$ and all
$K_i\in C^2_+$, we define the general mixed $L_{\phi}^*$-affine
surface area by
\begin{equation*} as^*(\phi_1, K_1;
\cdots ; \phi_n, K_n)=\int _{S^{n-1}} \prod _{i=1}^n
\left[\frac{\phi_i({f_{K_i}(u)h_{K_i}^{n+1}(u)})}{h_{K_i}(u)^n}\right]^{\frac{1}{n}}\,d\s
(u).\end{equation*}  Similarly, let $as^*_{\phi}(K_1,\cdots ,
K_n)=as^*(\phi, K_1; \cdots ; \phi, K_n),$
$as_{\phi}^*(K)=as_{\phi}^*(K, \cdots, K)$, and $as^*(\phi_1,
\cdots, \phi _n; K)=as^*(\phi_1, K; \cdots ; \phi_n, K)$. In
particular, $as^*(\phi_1, \cdots, \phi_n;
B^n_2)=\big[\phi_1(1)\cdots \phi_n(1)\big]^{\frac{1}{n}} n
|B^n_2|.$

\vskip 1mm  \begin{proposition}\label{dual:formula:star} Let
$K_i\in C^2_+$ be convex bodies, such that, $K_i=\l_{i}K$ for some
convex body $K\in C^2_+$ and $\l_i>0$, $i=1, \cdots, n$. For all
$\phi_i\in Conc(0,\infty)$,
$$as^*(\phi_1, K_1; \cdots ; \phi_n, K_n)=as(\phi_1,
K_1^\circ;\cdots ; \phi_n, K_n^\circ ).$$ In particular,
$as^*(\phi_1, \cdots, \phi_n; K)=as(\phi_1, \cdots, \phi_n;
K^\circ).$ \end{proposition}

 \vskip 1mm \noindent {\bf Remark.} The proof of this proposition is similar to that of
Proposition \ref{Duality:psi}. An immediate consequence is
$as_{\phi}^*(K)\!=\!as_{\phi}(K^\circ)$ for $K\!\in\! C^2_+$ and
$\phi\!\in\! Conc(0,\infty).$ In particular, if
$\phi(t)\!=\!t^{\frac{p}{n+p}}$ for $p\!\geq\!0$, one obtains the
duality formula $as_p(K)\!=\!as_{\frac{n^2}{p}}(K^\circ)$
\cite{Hug, WY}. Hence the $L_{\phi}^*$-affine surface area can be
viewed as a generalization of the $L_p$ affine surface area for
$p\!\geq\!0$. As above, one {\em cannot} expect, in general,
$as^*(\phi_1, K_1; \cdots ; \phi_n, K_n)=as(\phi_1,
K_1^\circ;\cdots ; \phi_n, K_n^\circ )$ even if all $\phi_i=\phi$
for some $\phi\in Conc(0, \infty)$ and all $K_i$ are ellipsoids.


\vskip 1mm The following theorem gives a geometric interpretation
for the general mixed $L_\phi^*$-affine surface area. Similar
geometric interpretations can also be obtained by the surface body
\cite{SW5, WY}. \bt Let $K, K_i\in C^2_+$ and $\psi_i\in
Conv(0,\infty)$, $i=1, \cdots, n$. Let $g:
\partial K\rightarrow \bbR$ be the function
\begin{equation*}
g(N_K^{-1}(u))=f_K(u)^{\frac{n-2}{2}}\prod_{i=1}^n\left[\frac{\psi_i({f_{K_i}(u)h_{K_i}^{n+1}(u)})}{h_{K_i}(u)^n}\right]^{\frac{1-n}{2n}}.
\end{equation*} Let $c_n=2|B^{n-1}_2|^{\frac{2}{n-1}}$.
 Then,
$$as^*(\psi_1, K_1; \cdots ; \psi_n, K_n)=\lim_{s\rightarrow 0} c_n\
\frac{|K^{g,s}|-|K|}{s^{\frac{2}{n-1}}}.$$ \et

\vskip 1mm \noindent {\bf Remark.} Replacing $\psi _i\in Conv(0,
\infty)$ by $\phi_i\in Conc(0, \infty)$, one gets the geometric
interpretation for the general mixed $L_{\phi}^*$-affine surface
area. That is,
$$as^*(\phi_1, K_1; \cdots ; \phi_n, K_n)=\lim_{s\rightarrow 0} c_n\
\frac{|K^{\tilde{g},s}|-|K|}{s^{\frac{2}{n-1}}},$$ where the
function $\tilde{g}$ takes the form
$$\tilde{g}(N_K^{-1}(u))=f_K(u)^{\frac{n-2}{2}}\prod_{i=1}^n\left[\frac{\phi_i({f_{K_i}(u)h_{K_i}^{n+1}(u)})}{h_{K_i}(u)^n}\right]^{\frac{1-n}{2n}}.$$

\vskip 1mm The general mixed $L_{\psi}^*$-affine surface area is
also affine invariant. For all $\psi_i=\psi$ and $K_i=K$, this
result was proved in \cite{Ludwig2009}. The proof is similar to
that of Theorem \ref{affine:invariant:1}, and we omit it.

\bt Let $T: \bbR^n \rightarrow\bbR^n $ be an invertible linear
transform. For all $\psi_i\in Conv(0,\infty)$ and $K_i\in C^2_+$,
one has
$$as^*(\psi_1, TK_1; \cdots ; \psi_n, TK_{n})=as^*(\psi_1, K_1; \cdots ; \psi_n,
K_{n}), \ for \ \ \ |det(T)|=1.$$ \noindent If in addition, all
$\psi_i \in Conv(0, \infty)$ are homogeneous of degrees $r_i\in
(-\infty, 0]$, then
$$as^*(\psi_1, TK_1; \cdots ; \psi_n, TK_{n})=|det(T)|^{\frac{2\sum
_{i=1}^n r_i}{n}-1} as^*(\psi_1, K_1; \cdots ; \psi_n, K_{n}).$$
\et

\noindent {\bf Remark.} Similar results hold for the general mixed
$L_{\phi}^*$-affine surface area, i.e.,
$$as^*(\phi_1, TK_1; \cdots ; \phi_n, TK_{n})=as^*(\phi_1, K_1; \cdots ; \phi_n,
K_{n}), \ for \ \ \ |det(T)|=1.$$ \noindent If in addition,
$\phi_i\in Conc(0, \infty)$ are homogeneous of degrees $r_i\in
[0,1)$, then
$$as^*(\phi_1, TK_1; \cdots ; \phi_n, TK_{n})=|det(T)|^{\frac{2\sum
_{i=1}^n r_i}{n}-1} as^*(\phi_1, K_1; \cdots ; \phi_n, K_{n}).$$


\section{Inequalities}


\vskip 1mm A general version of the classical Alexandrov-Fenchel
inequalities for  mixed volumes (see \cite{ Ale1937, Bus1958,
Sch}) can be written as
$$
\prod_{i=0}^{m-1}V(K_1,\cdots, K_{n-m}, \underbrace{K_{n-i},
\cdots, K_{n-i}}_m)\leq V^m(K_1, \cdots, K_n). $$ Here the
analogous inequalities for general mixed affine surface areas are
proved. We refer readers to the references \cite{Lut1975, Lut1987,
Lu1, WernerYe2010} for similar results related to the mixed
$p$-affine surfaces area.

\begin{proposition} \label{inequality:mixed:conc}
Let all $K_i\in C^2_+$ and $\phi_i\in Conc(0,\infty)$. Then, for
$1 \leq m \leq n$,
$$ as^m(\phi_1,\! K_1;\! \cdots\! ;\! \phi_n,\! K_{n})\! \leq \!\!\prod_{i=0}^
{m-1}as(\phi_1,\! K_1;\! \cdots;\! \phi_{n-m},\! K_{n-m};
\underbrace{\phi_{n-i},\! K_{n-i};\! \cdots;\! \phi_{n-i},\!
K_{n-i}}_m).
$$ Equality holds if (1) all $K_i$ coincide and $\phi_i=\l_i\phi_n$
for some $\l_i>0$, $i=n-m+1, \cdots, n$, or (2)
$\phi_i=\l_i\phi_n$, $K_i=\eta _i K_n$ for some $\l_i, \eta_i>0$,
$i=n-m+1, \cdots, n$, and $\phi_n$ is homogeneous of degree $r\in
[0,1)$. The equality holds trivially if $m=1$.

\vskip 1mm \noindent In particular, if $m=n$,
\begin{eqnarray}\label{inequality:mixed:conc:1} as^n(\phi_1,K_1;
\cdots; \phi_n, K_{n}) \leq  as_{\phi_1}(K_1) \cdots
as_{\phi_n}(K_n).
\end{eqnarray}
\end{proposition}

\vskip 1mm \noindent  {\bf Proof.} Let us put $$g_0(u)= \prod
_{i=1}^{n-m}
\!\left[\phi_i\bigg(\!\frac{1}{f_{K_i}(u)h_{K_i}^{n+1}(u)}\!\bigg)h_{K_i}(u)f_{K_i}(u)\right]^{\frac{1}{n}},$$
and for $j=0, \cdots, m-1$, put $$g_{j+1}(u)=
\!\left[\phi_{n-j}\bigg(\!\frac{1}{f_{K_{n-j}}(u)h_{K_{n-j}}^{n+1}(u)}\!\bigg)h_{K_{n-j}}(u)f_{K_{n-j}}(u)\right]^{\frac{1}{n}}.$$
By H\"{o}lder's inequality (see \cite{HLP})
\begin{eqnarray*}
as(\phi_1,\! K_1;\! \cdots\! ;\! \phi_n,\! K_n)\!\!&=&\!\!\!\!
\int
_{S^{n-1}}\! g_0(u) g_1(u) \cdots g_{m}(u)\,d\s(u)\\
\!\!\!\! &\leq&\!\! \prod _{j=0}^{m-1} \left(\int _{S^{n-1}}
g_0(u) g_{j+1}^m(u)\,d\s(u)\right)^{\frac{1}{m}}\\ \!\!&=&\!\!\!\!
\prod_{j=0}^ {m-1} as^{\frac{1}{m}} (\phi_1,\! K_1;\! \cdots\! ;
\!\phi_{n-m},\! K_{n-m};\! \underbrace{\phi_{n-j},\! K_{n-j};\!
\cdots;\! \phi_{n-j},\! K_{n-j}}_m).
\end{eqnarray*}

As $K_i \in C^2_+$ and $\phi_i\neq 0$, $g_k(u)>0$ for all $k=0, 1,
\cdots, m$ and all $u\in S^{n-1}$. Therefore, equality in
H\"{o}lder's inequality holds if and only if $g_0(u)
g_{j+1}^m(u)=\l_{n-j} ^m g_0(u) g_{1}^m(u)$ for some $\l_{n-j}>0$
and all $0\leq j\leq m-1$. This condition holds true if (1) all
$K_i$ coincide and $\phi_i=\l_i\phi_n$ for some $\l_i>0$,
$i=n-m+1, \cdots, n$, or (2) $\phi_i=\l_i\phi_n$, $K_i=\eta _i
K_n$ for some $\l_i, \eta_i>0$, $i=n-m+1, \cdots, n$, and $\phi_n$
is homogeneous of degree $r\in [0,1)$.

\vskip 1mm \noindent {\bf Remark.} When all
$\phi_i(t)=t^{\frac{p}{n+p}}$ for $p\geq 0$, this recovers the
Alexandrov-Fenchel inequalities for mixed $p$-affine surface areas
with $p\geq 0$ \cite{WernerYe2010}. When all $K_i$ coincide with
$K$, then $$as^m(\phi_1, \cdots, \phi_n; K)\leq \prod
_{i=0}^{m-1}as(\phi_1, \cdots, \phi_{n-m}, \underbrace{\phi_{n-i},
\cdots, \phi_{n-i}}_m; K),$$ and equality holds true if
$\phi_i=\l_i \phi_n$ for some $\l_i>0$, $i=n-m+1, \cdots, n$. When
all $\phi_i$ coincide with $\phi$, then
$$as^m_{\phi}(K_1, \cdots, K_n)\leq \prod
_{i=0}^{m-1}as_{\phi}(K_1, \cdots, K_{n-m}, \underbrace{K_{n-i},
\cdots, K_{n-i}}_m),$$ and equality holds true if (1) all $K_i$,
$i=n-m+1, \cdots, n$ coincide, or (2) $K_i=\eta _i K_n$ for some
$\eta _i>0$, $i=n-m+1, \cdots, n$ and $\phi$ is homogeneous of
degree $r\in [0,1)$.

\vskip 1mm \noindent {\bf Remark.}  With slight modification, one
can get analogous Alexandrov-Fenchel inequality for the general
mixed $L^*_{\phi}$-affine surface area, namely,
$$ \big[as^*(\phi_1,\! K_1;\! \cdots\! ; \phi_n,\! K_{n})\big]^m
\!\!\leq\!\! \prod_{i=0}^ {m-1}\! as^*(\phi_1,\! K_1;\! \cdots\! ;
\phi_{n-m},\! K_{n-m};\! \underbrace{\phi_{n-i},\! K_{n-i};\!
\cdots;\! \phi_{n-i},\! K_{n-i}}_m). $$ In particular, if $m=n$,
\begin{eqnarray}\label{inequality:mixed:convex:1} \big[as^*(\phi_1,K_1;
\cdots; \phi_n, K_{n})\big]^n \leq  as_{\phi_1}^*(K_1) \cdots
as_{\phi_n}^*(K_n).
\end{eqnarray} Replacing $\phi_i\in
Conc(0,\infty)$ by $\psi_i\in Conv(0,\infty)$, one gets the
Alexandrov-Fenchel inequalities for the general mixed $L_{\psi}$-
and $L^*_{\psi}$-affine surface areas.


\vskip 1mm Blaschke-Santal\'{o} inequality states that, for all
convex body $K$ with centroid at the origin, $|K||K^\circ| \leq
|B^n_2|^2$; equality holds if and only if $K$ is an ellipsoid.
This fundamental inequality has been generalized to the $L_p$
affine surface area and mixed $p$-affine surface area \cite{Lu1,
WY, WernerYe2010}. Here, we prove the Santal\'{o}-type
inequalities for the general mixed $L_{\phi}$-affine surface area.
The Santal\'{o}-type inequality for the general mixed
$L_{\phi}^*$-affine surface area can be achieved in the same way.

\bt\label{Stantalo:phi} Let $K_i\in C^2_+$ and $\phi_i \in
Conc(0,\infty)$ be homogeneous of degrees $r_i\in [0,1)$. Then
$$as(\phi_1, K_1; \cdots; \phi_n, K_n)as(\phi_1, K_1^\circ; \cdots
; \phi_n, K_n^\circ)\leq n^2 \prod_{i=1}^n \big[\phi_i(1)^2\
|K_i|\ |K_i^\circ|\big]^{\frac{1}{n}}.$$ Moreover, if all $K_i\in
C^2_+$ have centroid at the origin, then
$$as(\phi_1, K_1; \cdots; \phi_n, K_n)as(\phi_1, K_1^\circ;
\cdots ; \phi_n, K_n^\circ)\leq \big[as(\phi_1, B^n_2; \cdots;
\phi_n, B^n_2)\big]^2,$$ with equality if all $K_i$ are
origin-symmetric ellipsoids that are dilates of one another. \et

 \vskip 1mm \noindent {\bf
Proof.} The following inequality was proved in \cite{Ludwig2009}
and it actually holds true for all $K\in C^2_+$: if $\phi \in
Conc(0,\infty)$, then
\begin{equation} \label{Monika:inequality}
as_{\phi}(K)\leq n|K| \phi \bigg(\frac{|K^\circ|}{|K|}\bigg).
\end{equation}
Thus, for any $K\in C^2_+$,
\begin{eqnarray*}  as_{\phi}(K)as_{\phi}(K^\circ)\leq  n^2|K||K^\circ| \phi
\bigg(\frac{|K^\circ|}{|K|}\bigg)\phi
\bigg(\frac{|K|}{|K^\circ|}\bigg).\end{eqnarray*}

\vskip 1mm \noindent Combining with inequality
(\ref{inequality:mixed:conc:1}), one has, for all $K_i\in C^2_+$,
\begin{eqnarray} as(\phi_1,\! K_1;\! \cdots;\! \phi_n,\!
K_n)as(\phi_1,\! K_1^\circ;\! \cdots ;\! \phi_n,\! K_n^\circ) \leq
\prod_{i=1}^n \big[as_{\phi_i}(K_i)as_{\phi _i}(K^\circ_i)\big]^{\frac{1}{n}}\nonumber \\
\leq n^2\! \prod_{i=1}^n\! \left[\!|K_i| |K_i^\circ|\phi_i\!
\bigg(\!\frac{|K_i^\circ|}{|K_i|}\!\bigg)\phi_i \!
\bigg(\!\frac{|K_i|}{|K_i^\circ|}\!\bigg)\right]^{\frac{1}{n}}=n^2\!
\prod_{i=1}^n \big[\phi_i(1)^2\ |K_i|\
|K_i^\circ|\big]^{\frac{1}{n}}. \label{santalo:conc:1}
\end{eqnarray} Here the equality follows from all $\phi_i$
being of homogenous degrees $r_i\in [0,1)$.

\vskip 1mm If all $K_i\in C^2_+$ have centroid at the origin, one
can employ Blaschke-Santal\'{o} inequality to inequality
(\ref{santalo:conc:1}) and get \begin{eqnarray*} as(\phi_1,\!
K_1;\! \cdots\!; \phi_n,\! K_n)as(\phi_1,\! K_1^\circ;\! \cdots\!
; \phi_n,\! K_n^\circ)\!\!\leq\!\! n^2\! \prod_{i=1}^n
\big[\phi_i(1)\!
|B^n_2|\big]^{\frac{2}{n}}\!\!=\!\!\big[\!as(\phi_1,\! \cdots\!,
\phi_n\!; B^n_2)\!\big]^2.
\end{eqnarray*}
Clearly, equality holds true if $K_i$ are origin-symmetric
ellipsoids that are dilates of one another.


\vskip 1mm  Let $K\in C^2_+$ be a convex body with centroid at the
origin and $B_K$ be the origin-symmetric (Euclidean) ball such
that $|B_K|=|K|$. For the $L_{\phi}$-affine surface area with
$\phi \in Conc(0,\infty)$, Ludwig proved the affine isoperimetric
inequality \cite{Ludwig2009}; namely, $as_{\phi}(K)\leq
as_{\phi}(B_K)$ with equality if and only if $K$ is an ellipsoid.
If we further assume that $\phi$ is homogeneous of degree $r\in
[0, 1)$, then the affine isoperimetric inequality for
$L_{\phi}$-affine surface area may be stated as
\begin{eqnarray}\label{affine:isoperimetric:inequality:homo}\left(\frac{as_{\phi}(K)}{as_{\phi}(B^n_2)}\right)\leq
\left(\frac{|K|}{|B^n_2|}\right)^{1-2r},\end{eqnarray} with
equality if and only if $K$ is an ellipsoid. In fact, let
$\l=\frac{|B^n_2|}{|K|}$ and $\tilde{K}=\l^{\frac{1}{n}} K$, then
$|\tilde{K}|=|B^n_2|$. Employing the affine isoperimetric
inequality in \cite{Ludwig2009} to $\tilde{K}$, one has
$\frac{as_{\phi}(\tilde{K})}{as_{\phi}({B^n_2})}\leq 1.$ By
Theorem \ref{affine:invariant:1}, one has
$$\frac{as_{\phi}(\tilde{K})}{as_{\phi}({B^n_2})}=\frac{as_{\phi}(\l^{\frac{1}{n}}{K})}{as_{\phi}({B^n_2})}=
\l ^{1-2r}\left(\frac{as_{\phi}({K})}{as_{\phi}({B^n_2})}\right)
\leq 1,
$$ which is equivalent to the formula (\ref{affine:isoperimetric:inequality:homo}). There is an
equality if and only if $K$ is an ellipsoid. When
$\phi(t)=t^{\frac{p}{n+p}}$ for $p>0$, one gets the $L_p$ affine
isoperimetric inequality for $p>0$ \cite{Lu1, WY}. For $p=0$, one
has equality instead of inequality.


 \vskip 1mm Next, we prove the affine isoperimetric
inequalities for general mixed affine surface areas. Hereafter, we
always let $B_{K_i}$ be the origin-symmetric (Euclidean) ball s.t.
$|B_{K_i}|=|K_i|$ for all $i$.

\bt \label{affine:isoperimetric:inequality} Let all $K_i\in C^2_+$
be convex bodies with centroid at the origin. \vskip 1mm \noindent
(i): If all $\phi_i\in Conc(0, \infty)$, then
$${as(\phi _1, K_1; \cdots ; \phi_n, K_n)}\leq {as(\phi_1, B_{K_1}; \cdots; \phi_n, B_{K_n})}.$$
Equality holds if all $K_i$ are ellipsoids that are dilates of one
another.

\vskip 1mm \noindent (ii) For all $\phi_i\in Conc(0,\infty)$ with
homogeneous degrees $r_i\in [0,1)$,
$$\left(\frac{as(\phi _1, K_1; \cdots ; \phi_n, K_n)}{as(\phi_1,
\cdots, \phi_n; B^n_2)}\right)^n \leq \prod _{i=1}^n
\left(\frac{|K_i|}{|B^n_2|}\right)^{1-2r_i},$$ with equality if
all $K_i$ are ellipsoids that are dilates of one another. \et

\vskip 1mm \noindent {\bf Proof.}

\vskip 1mm \noindent (i). Let all $K_i\!\in\! C^2_+$ be convex
bodies with centroid at the origin and $|K_i|\!=\!|B_{K_i}|$ for
all $i$. It is easy to verify that $as_{\phi_1}(B_{K_1})\!
\cdots\! as_{\phi_n}(B_{K_n})\!=\!\big[as(\phi_1,\! B_{K_1};\!
\cdots;\! \phi_n,\! B_{K_n})\big]^n$. By inequality
(\ref{inequality:mixed:conc:1}) and the affine isoperimetric
inequality for $L_{\phi}$-affine surface area, one has, for all
$\phi_i\in Conc(0,\infty)$,
\begin{eqnarray*}as^n(\phi_1,\!K_1;\!
\cdots;\! \phi_n,\! K_{n})\!&\leq&\! as_{\phi_1}(K_1) \cdots
as_{\phi_n}(K_n)\!\leq\! as_{\phi_1}(B_{K_1})\! \cdots\!
as_{\phi_n}(B_{K_n})\\ \!&=&\!\big[as(\phi_1,\! B_{K_1};\!
\cdots;\! \phi_n,\! B_{K_n})\big]^n.\end{eqnarray*} By the affine
invariant property in Theorem \ref{affine:invariant:1}, the
equality holds true if $K_i$ are ellipsoids that are dilates of
one another.

\vskip 1mm \noindent (ii). Recall that $\big[as(\phi_1,\!
\cdots\!, \phi_n;\! B^n_2)\big]^n\!=\! as_{\phi_1}(B_2^n)\!
\cdots\! as_{\phi_n}(B_2^n)$. By inequalities
(\ref{inequality:mixed:conc:1}) and
(\ref{affine:isoperimetric:inequality:homo}), one has
\begin{eqnarray*}\left(\frac{as(\phi_1,\! K_1;\!
\cdots\!; \phi_n,\! K_{n})}{as(\phi_1,\! \cdots\!, \phi_n;\!
B^n_2)}\right)^n \!&\leq &\! \frac{as_{\phi_1}(K_1)\! \cdots \!
as_{\phi_n}(K_n)}{as_{\phi_1}(B_2^n)\! \cdots\!
as_{\phi_n}(B_2^n)}\!\leq\! \prod _{i=1}^n
\left(\frac{|K_i|}{|B^n_2|}\right)^{1-2r_i}.
\end{eqnarray*} If all $K_i$ are ellipsoids that are dilates of one
another, the equality holds true.

\vskip 1mm \noindent \begin{remark}\label{remark:1} If all
$\phi_i(t)=t^{\frac{p}{n+p}}$ for $p\geq 0$, one recovers affine
isoperimetric inequalities for mixed $p$-affine surface areas
\cite{WernerYe2010}. One cannot expect to get strictly positive
lower bounds in Theorem \ref{affine:isoperimetric:inequality}. Let
the convex body $K(R, \varepsilon)\subset \bbR^2$ be the
intersection of four Euclidean balls with radius $R$ centered at
$(\pm (R-1), 0)$, $(0, \pm (R-1))$, $R$ arbitrarily large. We then
``round" the corners by putting there arcs of Euclidean balls of
arbitrarily small radius $\varepsilon$, and ``bridge" between the
$R$-arcs and $\varepsilon$-arcs by $C^2_+$-arcs on a set of
arbitrarily small measure to obtain a convex body in $C^2_+$
\cite{WY}. Then $as_{\phi}(K(R, \varepsilon))\leq 16
R^{-\frac{p}{p+2}}+4\pi \varepsilon ^{\frac{2}{2+p}}$ for
$\phi(t)=t^{\frac{p}{n+p}}$ with $p>0$, which goes to $0$ as
$R\rightarrow \infty$ and $\varepsilon \rightarrow 0$. Choose now
$R_i$ and $\varepsilon_i$, $1\leq i\leq n$, such that
$R_i\rightarrow \infty$ and $\varepsilon_i\rightarrow 0$, and let
$K_i=K(R_i, \varepsilon_i)$ for $i=1, \cdots, n$. Let
$\phi_i(t)=t^{\frac{p_i}{n+p_i}}$ for $p_i>0$, by inequality
(\ref{inequality:mixed:conc:1}), $as^n(\phi_1, K_1; \cdots ;
\phi_n, K_n)\leq \prod _{i=1}^n as_{\phi_i}(K_i),$ and thus
$as(\phi_1, K_1; \cdots ; \phi_n, K_n)\rightarrow 0$.\end{remark}

\vskip 1mm We can prove the following affine isoperimetric
inequality for the general mixed $L_{\phi}^*$-affine surface area.

\bt \label{affine:isoperimetric:inequality:star} Let all $K_i\in
C^2_+$ be convex bodies with centroid at the origin. \vskip 1mm
\noindent (i): If all $\phi_i\!\in\! Conc(0, \infty)$, then
$${as^*(\phi _1,\! K_1;\! \cdots ;\! \phi_n,\! K_n)}\!\leq\!
as^*(\phi_1,\! (B_{K_1^\circ})^\circ;\! \cdots;\! \phi_n,\!
(B_{K_n^\circ})^\circ).$$ Equality holds if $K_i$ are ellipsoids
that are dilates of one another.

\vskip 1mm \noindent (ii) For all $\phi_i\in Conc(0,\infty)$ with
homogeneous degrees $r_i\in [0, 1)$,
$$\left(\frac{as^*(\phi _1, K_1; \cdots ; \phi_n, K_n)}{as^*(\phi_1,
\cdots, \phi_n; B^n_2)}\right)^n \leq \prod _{i=1}^n
\left(\frac{|K_i|}{|B^n_2|}\right)^{2r_i-1},$$ with equality if
$K_i$ are ellipsoids that are dilates of one another.   \et

\vskip 1mm \noindent {\bf Proof.} \vskip 1mm \noindent (i). By
inequality (\ref{inequality:mixed:convex:1}), Proposition
\ref{dual:formula:star}, and the affine isoperimetric inequality
for the $L_{\phi}$-affine surface area, one gets,
\begin{eqnarray*}
[as^*(\phi_1,K_1; \cdots; \phi_n, K_{n})]^n \!\!&\leq &\!\!
as_{\phi_1}^*(K_1) \cdots as_{\phi_n}^*(K_n) = as_{\phi_1}(K_1^\circ) \cdots as_{\phi_n}(K_n^\circ)\\
\!\!&\leq&\!\! as_{\phi_1}(B_{K_1^\circ}) \cdots
as_{\phi_n}(B_{K_n^\circ}) \!\!=\!\! \big[as(\phi_1,\!
B_{K_1^\circ};\! \cdots;\! \phi_n,\! B_{K_n^\circ})\big]^n\\
\!\!&=&\!\! \big[as^*(\phi_1,\! (B_{K_1^\circ})^\circ;\! \cdots;\!
\phi_n,\! (B_{K_n^\circ})^\circ)\big]^n,
\end{eqnarray*} where the last equality follows
Proposition \ref{dual:formula:star}. Following the affine
invariant property, equality holds if $K_i$ are ellipsoids that
are dilates of one another.

\vskip 1mm \noindent (ii). Similarly, by inequalities
(\ref{inequality:mixed:convex:1}), (\ref{Monika:inequality}), and
Proposition \ref{dual:formula:star}, one has
\begin{eqnarray*}\big[as^*(\phi_1,\! K_1;\! \cdots\!; \phi_n,\!
K_{n})\big]^n \leq {as_{\phi_1}^*(K_1) \cdots as_{\phi_n}^*(K_n)}
= {as_{\phi_1}(K_1^\circ) \cdots as_{\phi_n}(K_n^\circ)}\\
\leq n^n \prod _{i=1}^n \big[|K_i^\circ|
\phi_i(|K_i|/|K^\circ_i|)\big] = n^n \prod _{i=1}^n
\big[\phi_i(1)|K_i^\circ|^{1-r_i}|K_i|^{r_i} \big].\nonumber
\end{eqnarray*} By Blaschke-Santal\'{o}
inequality, i.e., $|K^\circ||K|\leq  |B^n_2|^2$, and $r_i\in
[0,1)$, one gets
\begin{eqnarray*}\big[as^*(\phi_1,\! K_1;\! \cdots\!; \phi_n,\!
K_{n})\big]^n \leq n^n \prod_{i=1}^n
\big[\phi_i(1)|K_i|^{2r_i-1}|B^n_2|^{2-2r_i}\big].\end{eqnarray*}
Equivalently, by $as^*(\phi_1, \cdots, \phi_n;
B^n_2)=\big[\phi_1(1)\cdots \phi_n(1)\big]^{\frac{1}{n}} n
|B^n_2|,$ one has
\begin{eqnarray*}\left(\frac{as^*(\phi_1,\! K_1;\! \cdots\!; \phi_n,\!
K_{n})}{as^*(\phi_1, \cdots, \phi_n; B^{n}_2)}\right)^n \leq
\prod_{i=1}^n
\left(\frac{|K_i|}{|B^n_2|}\right)^{2r_i-1}.\end{eqnarray*}
Clearly, equality holds true if all $K_i$ are ellipsoids that are
dilates of one another.

\vskip 1mm \noindent \begin{remark} One cannot expect to get
strictly positive lower bounds in Theorem
\ref{affine:isoperimetric:inequality:star}. Let the convex bodies
$K_i=K(R_i, \varepsilon_i)\in C^2_+$ be as in Remark
(\ref{remark:1}) with $R_i\rightarrow \infty$ and
$\varepsilon_i\rightarrow 0$. Let
$\phi_i(t)\!=\!t^{\frac{p_i}{n+p_i}}$ with $p_i\!>\!0, i\!=\!1,
\!\cdots,\! n$. By inequality (\ref{inequality:mixed:convex:1})
and Proposition \ref{Duality:psi}, one gets $\big[as^{*}(\phi_1,\!
K_1^\circ;\! \cdots;\! \phi_n,\! K_n^\circ)\big]^n\!\leq\!
\prod_{i=1}^n as_{\phi_i}^*(K_i^\circ)\!=\!\prod_{i=1}^n
as_{\phi_i}(K_i),$ which goes to $0$ as $R_i\rightarrow \infty$
and $\varepsilon\rightarrow 0 $. \end{remark}


\section{General $i$-th mixed affine surface areas and related
inequalities} Let $i$ be a real number and $\phi_1, \phi_2\in
Conc(0, \infty)$. The general $i$-th mixed $L_{\phi}$-affine
surface area of $K, L \in C^2_+$ is defined as
\begin{equation}\label{i:mixed:phi} as_{i}\!(\phi_1,\! K;\! \phi_2,\!
L)\!\!=\!\!\!\int
_{S^{n-1}}\!\!\left[\!\phi_1\!\bigg(\!\frac{f^{-1}_{K}(u)}{h_{K}^{n+1}(u)}\!\bigg)h_{K}(u)f_{K}(u)\!\right]^{\frac{n-i}{n}}\!\!
\left[\!\phi_2\bigg(\!\frac{f^{-1}_{L}(u)}{h_{L}^{n+1}(u)}\!\bigg)\!h_{L}(u)\!f_{L}(u)\!\right]^{\frac{i}{n}}\!\!\,d\s(u).
\end{equation} If $i$ is an integer number with $0\leq i\leq n$, then $$as_i(\phi_1, K; \phi_2, L)=as(\underbrace{\phi_1, K;\cdots ;
\phi_1, K}_{n-i}; \underbrace{\phi_2, L; \cdots; \phi_2, L}_i).$$
In particular, $as_i(\phi_1, K; \phi_2, L)=as_{n-i}(\phi_2, L;
\phi_1, K)$, $as_0(\phi_1, K; \phi_2, L)=as_{\phi_1}(K)$, and
$as_n(\phi_1, K; \phi_2, L)=as_{\phi_2}(L)$. If
$\phi_1(t)=\phi_2(t)=t^{\frac{p}{n+p}}$ for $p\geq 0$, one gets
the $i$-th mixed $p$-affine surface area \cite{Lut1987, Wa2007,
WernerYe2010}. Note that the $i$-th mixed $p$-affine surface area
includes the surface area of $K$ as a special case; namely the
surface area of $K$ is $as_{-1}(\phi, K; \phi, B^n_2)$ with
$\phi(t)=t^{\frac{1}{n+1}}$. Obviously, the general $i$-th mixed
$L_{\phi}$-affine surface area satisfies the affine invariant
property as in Theorem \ref{affine:invariant:1}. Similarly, we can
define the general $i$-th mixed $L_{\phi}^*$-affine surface area
of $K, L\in C^2_+$ as $$ as_{i}^*\!(\phi_1,\! K;\! \phi_2,\!
L)\!\!=\!\!\!\int
_{S^{n-1}}\!\!\left[\frac{\phi_1({f_{K}(u)h_{K}^{n+1}(u)})}{h_{K}^n(u)}\right]^{\frac{n-i}{n}}\!
\left[\frac{\phi_2({f_{L}(u)h_{L}^{n+1}(u)})}{h_{L}^n(u)}\right]^{\frac{i}{n}}\,d\s(u).
$$

\vskip 1mm For $\psi_1, \psi_2\in Conv(0,\infty)$, the general
$i$-th mixed $L_{\psi}$-affine surface area of $K, L \in C^2_+$ is
defined as
$$ as_{i}\!(\psi_1,\! K;\! \psi_2,\! L)\!\!=\!\!\!\int
_{S^{n-1}}\!\!\left[\!\psi_1\!\bigg(\!\frac{f^{-1}_{K}(u)}{h_{K}^{n+1}(u)}\!\bigg)h_{K}(u)f_{K}(u)\!\right]^{\frac{n-i}{n}}\!\!
\left[\!\psi_2\bigg(\!\frac{f^{-1}_{L}(u)}{h_{L}^{n+1}(u)}\!\bigg)\!h_{L}(u)\!f_{L}(u)\!\right]^{\frac{i}{n}}\!\!\,d\s(u),
$$ and the general $i$-th mixed $L_{\psi}^*$-affine surface area
of $K, L\in C^2_+$ is defined as
$$as_{i}^*\!(\psi_1,\! K;\! \psi_2,\! L)\!\!=\!\!\!\int
_{S^{n-1}}\!\!\left[\frac{\psi_1({f_{K}(u)h_{K}^{n+1}(u)})}{h_{K}^n(u)}\right]^{\frac{n-i}{n}}\!
\left[\frac{\psi_2({f_{L}(u)h_{L}^{n+1}(u)})}{h_{L}^n(u)}\right]^{\frac{i}{n}}\,d\s(u).
$$ If $\psi_1(t)=\psi_2(t)=t^{\frac{p}{n+p}}$ for $-n<p\leq 0$, then
$as_{i}(\psi_1,\! K;\! \psi_2,\! L)$ equals to the $i$-th mixed
$p$-affine surface area for $-n<p\leq 0$. On the other hand, if
$\psi_1(t)=\psi_2(t)=t^{\frac{n}{n+p}}$ for $p<-n$, then
$as_{i}^*(\psi_1,\! K;\! \psi_2,\! L)$ equals to the $i$-th mixed
$p$-affine surface area for $p<-n$ \cite{WernerYe2010}.

\begin{proposition}
If $j<i<k$ or $k<i<j$ (equivalently, $\frac{k-j}{k-i}>1$), then
for all $\phi_1, \phi_2\in Conc(0, \infty)$ and $K, L\in C^2_+$,
$$as_{i}(\phi_1, K; \phi_2, L)\leq as_{j}(\phi_1, K; \phi_2,
L)^{\frac{k-i}{k-j}}as_{k}(\phi_1, K; \phi_2,
L)^{\frac{i-j}{k-j}}.$$ Equality holds if (1) $K\!=\!L$ and
$\phi_1\!=\!\l \phi_2$ for some $\l\!>\!0$; (2) $\phi_1$ is
homogeneous of degree $r\in [0,1)$, $\phi_1\!=\!\l \phi_2$ for
some $\l\!>\!0$, and $K$ and $L$ are dilates of each other.
\end{proposition}

\noindent {\bf Remark.} Similar results can be obtained for other
general $i$-th mixed affine surface areas, for instance,
$as_{i}(\psi_1, K; \psi_2, L)\leq as_{j}(\psi_1, K; \psi_2,
L)^{\frac{k-i}{k-j}}as_{k}(\psi_1, K; \psi_2,
L)^{\frac{i-j}{k-j}}.$

\vskip 1mm \noindent {\bf Proof.} By formula (\ref{i:mixed:phi}),
one has \begin{eqnarray*} as_{i}\!(\phi_1,\! K;\! \phi_2,\! L)
\!\!\!\!&=&\!\!\!\! \int
_{S^{n-1}}\!\!\left[\!\phi_1\!\bigg(\!\frac{f^{-1}_{K}(u)}{h_{K}^{n+1}(u)}\!\bigg)h_{K}(u)f_{K}(u)\!\right]^{\frac{n-i}{n}}\!\!
\left[\!\phi_2\bigg(\!\frac{f^{-1}_{L}(u)}{h_{L}^{n+1}(u)}\!\bigg)\!h_{L}(u)\!f_{L}(u)\!\right]^{\frac{i}{n}}\!\!\,d\s(u)\\
\!\!\!\! &=&\!\!\!\!\int
_{S^{n-1}}\!\!\left\{\!\left[\!\phi_1\!\bigg(\!\frac{f^{-1}_{K}(u)}{h_{K}^{n+1}(u)}\!\bigg)h_{K}(u)f_{K}(u)\!\right]^{\frac{n-j}{n}}\!\!\!
\left[\!\phi_2\!\bigg(\!\frac{f^{-1}_{L}(u)}{h_{L}^{n+1}(u)}\!\bigg)\!h_{L}(u)\!f_{L}(u)\!\right]^{\frac{j}{n}}\!\right\}^{\frac{k-i}{k-j}}\\
\!\!\!\!&&\!\!\!\!\times\!
\!\left\{\!\left[\!\phi_1\!\bigg(\!\frac{f^{-1}_{K}(u)}{h_{K}^{n+1}(u)}\!\bigg)h_{K}(u)f_{K}(u)\!\right]^{\frac{n-k}{n}}\!\!\!
\left[\!\phi_2\!\bigg(\!\frac{f^{-1}_{L}(u)}{h_{L}^{n+1}(u)}\!\bigg)\!h_{L}(u)\!f_{L}(u)\!\right]^{\frac{k}{n}}\!\right\}^{\frac{i-j}{k-j}}
\!\!\!\!\,d\s(u) \\ \!\!\!\!&\leq& \!\!\!\! as_{j}(\phi_1, K;
\phi_2, L)^{\frac{k-i}{k-j}}as_{k}(\phi_1, K; \phi_2,
L)^{\frac{i-j}{k-j}},
\end{eqnarray*}
where the last inequality follows H\"{o}lder inequality and
formula (\ref{i:mixed:phi}). Clearly, the equality holds true if
(1) $K=L$ and $\phi_1=\l \phi_2$ for some $\l>0$; (2) $\phi_1$ is
homogeneous of degree $r\in [0,1)$, $\phi_1=\l \phi_2$ for some
$\l>0$, and $K, L$ are dilates of each other.

\vskip 1mm If $j=0, k=n$, then for all $0\leq i\leq n$,
\begin{equation}\label{i:mixed:phi:1}\big[as_i(\phi_1, K; \phi_2, L)\big]^n\leq
[as_{\phi_1}(K)]^{n-i}[as_{\phi_2}(L)]^{i}.\end{equation} On the
other hand, if $i=0$, $j=n$, $k\leq 0$, or $i=n, j=0$, $k\geq n$,
one has \begin{equation}\label{i:mixed:phi:2}\big[as_k(\phi_1, K;
\phi_2, L)\big]^n\geq
[as_{\phi_1}(K)]^{n-k}[as_{\phi_2}(L)]^{k}.\end{equation}

The following proposition gives the Santal\'{o}-type inequality
for the general $i$-th mixed $L_{\phi}$-affine surface area.
Similar results for the general $i$-th mixed $L_{\phi}^*$-affine
surface area also hold.

\begin{proposition} Let $0\!\leq\! i\!\leq\! n$, and $K,\! L\!\in\!
C^2_+$ be convex bodies with centroid at the origin. For
$\phi_1,\! \phi_2\!\in\! Conc(0, \infty)$ with homogeneous degrees
$r_1,\! r_2\!\in\! [0,1)$ respectively,
$$as_i(\phi_1, K; \phi_2, L)as_i(\phi_1, K^\circ; \phi_2,
L^\circ)\leq [as_i(\phi_1, B^n_2; \phi_2, B^n_2)]^2.$$ Equality
holds true if $K$ and $L$ are ellipsoids that are dilates of one
another.
\end{proposition}
\noindent {\bf Proof.} By inequalities (\ref{Monika:inequality})
and (\ref{i:mixed:phi:1}), one has
\begin{eqnarray*}as_i(\phi_1,\! K;\! \phi_2,\! L) as_i(\phi_1,\! K^\circ;\!
\phi_2,\! L^\circ )\!\! &\leq&\!\!
\big[as_{\phi_1}(K)as_{\phi_1}(K^\circ)]^{\frac{n-i}{n}}[as_{\phi_2}(L)as_{\phi_2}(L^\circ)\big]^{\frac{i}{n}}\\
\!\!&\leq&\!\!
n^2\big[|K||K^\circ|\phi_1(1)^2]^{\frac{n-i}{n}}[|L||L^\circ|\phi_2(1)^2\big]^{\frac{i}{n}}\\
\!\!&\leq&\!\! n^2
\phi_1(1)^{\frac{2(n-i)}{n}}\phi_2(1)^{\frac{2i}{n}} |B^n_2|^2\\
\!\!&=&\!\! \big[as_i(\phi_1, B^n_2; \phi_2, B^n_2)\big]^2,
\end{eqnarray*}
 where the last inequality follows from Blaschke-Santal\'{o} inequality and
 $0\leq i\leq n$. Clearly, the equality holds true if $K,L$ are
ellipsoids that are dilates of one another.

 \vskip 1mm The following proposition states the affine isoperimetric
 inequality for the general $i$-th mixed $L_{\phi}$-affine surface
 area.

\begin{proposition} Let $0\leq i\leq n$ and $K, L\in C^2_+$ be
convex bodies with centroid at the origin. For $\phi_1, \phi_2
\!\in\! Conc(0, \infty)$, one has

\noindent (i) $as_i(\phi_1, K; \phi_2, L) \leq as_i(\phi_1, B_K;
\phi_2, B_L)$, with equality if $K$ and $L$ are ellipsoids that
are dilates of each another;

\noindent (ii) if in addition, $\phi_1, \phi_2$ are homogeneous of
degrees $r_1, r_2\!\in\! [0,1)$ respectively,
$$\left(\frac{as_i(\phi_1, K; \phi_2, L)}{as_i(\phi_1, B^n_2; \phi_2, B^n_2)}\right)^n \leq
\left(\frac{|K|}{|B^n_2|}\right)^{(n-i)(1-2r_1)}\
\left(\frac{|L|}{|B^n_2|}\right)^{i(1-2r_2)}.$$ Equality holds
true if $K$ and $L$ are ellipsoids that are dilates of one
another.
\end{proposition}

\noindent {\bf Remark.} Similarly, one can get the affine
isoperimetric inequality for the general $i$-th mixed
$L^*_{\phi}$-affine surface area. For instance, if $\phi_1, \phi_2
\in Conc(0, \infty)$ with homogeneous degrees $r_1, r_2\in [0,1)$
respectively, then
$$\left(\frac{as_i^*(\phi_1, K; \phi_2, L)}{as_i^*(\phi_1, B^n_2; \phi_2, B^n_2)}\right)^n \leq
\left(\frac{|K|}{|B^n_2|}\right)^{(n-i)(2r_1-1)}\
\left(\frac{|L|}{|B^n_2|}\right)^{i(2r_2-1)}.$$ Equality holds if
$K$ and $L$ are ellipsoids that are dilates of one another.

\vskip 1mm \noindent {\bf Proof.} (i) Note $[as_i(\phi_1, B_K;
\phi_2, B_L)]^n=[as_{\phi_1}(B_K)]^{n-i}[as_{\phi_2}(B_L)]^i$. The
desired result is then an immediate consequence of inequality
(\ref{i:mixed:phi:1}), and Ludwig's isoperimetic inequality
\cite{Ludwig2009}.

\vskip 1mm\noindent (ii) By inequality (\ref{i:mixed:phi:1}), and
$[as_i(\phi_1, B^n_2; \phi_2, B^n_2)]^n=[as_{\phi_1}(B^n_2)]^{n-i}
[as_{\phi_2}(B^n_2)]^i$,
\begin{eqnarray*}\left[\frac{as_i(\phi_1, K; \phi_2,
L)}{as_i(\phi_1, B^n_2; \phi_2, B^n_2)}\right]^n &\leq&
\left[\frac{as_{\phi_1}(K)}{as_{\phi_1}(B^n_2)}\right]^{n-i}\left[\frac{as_{\phi_2}(L)}{as_{\phi_2}(B^n_2)}\right]^{i}.
\end{eqnarray*} Combining with inequality
(\ref{affine:isoperimetric:inequality:homo}) and $r_1, r_2\in
[0,1)$, the desired result follows. Clearly, equality holds if $K$
and $L$ are ellipsoids that are dilates of one another.

\begin{proposition}\label{ith:phi:bigger:n} Let $K\in C^2_+$ be a convex body with
centroid at the origin. For $k\geq n$, and $\phi_1, \phi_2 \in
Conc(0, \infty)$, one has

\vskip 1mm \noindent (i) $as_k(\phi_1,\! K;\! \phi_2,\!
B^n_2)\!\geq\! as_k(\phi_1,\! B_K;\! \phi_2,\! B^n_2),$ with
equality if $K$ is a ball;

\vskip 1mm \noindent (ii) if in addition, $\phi_1$ is homogeneous
of degree $r_1\in [0, 1)$, then
$$\left(\frac{as_k(\phi_1, K; \phi_2, B^n_2)}{as_k(\phi_1, B^n_2;
\phi_2, B^n_2)}\right)^n\geq
\left(\frac{|K|}{|B^n_2|}\right)^{(n-k)(1-2r_1)},$$ with equality
if $K$ is a ball; Moreover, $$as_k(\phi_1, K; \phi_2,
B^n_2)as_k(\phi_1, K^\circ; \phi_2, B^n_2)\geq [as_k(\phi_1,
B^n_2; \phi_2, B^n_2)]^2,$$ with equality if $K$ is a
ball.\end{proposition}

\vskip 1mm \noindent {\bf Proof.} (i) By inequality
(\ref{i:mixed:phi:2}) and $[as_k(\phi_1,\! B_K;\! \phi_2,\!
B^n_2)]^n\!=\![as_{\phi_1}(B_K)]^{n-k} [as_{\phi_2}(B^n_2)]^k$,
\begin{eqnarray*}\left[\frac{as_k(\phi_1,\! K;\! \phi_2,\!
B^n_2)}{as_k(\phi_1,\! B_K;\! \phi_2,\! B^n_2)}\right]^n &\geq&
\left[\frac{as_{\phi_1}(K)}{as_{\phi_1}(B_K)}\right]^{n-k}\geq 1
,\end{eqnarray*} where we have used the Ludwig's isoperimetric
inequality in \cite{Ludwig2009} and $n-k\leq 0$. The equality
holds trivially if $K$ is a ball.

\vskip 1mm \noindent (ii) Again, by inequality
(\ref{i:mixed:phi:2}) and $[as_k(\phi_1, B^n_2; \phi_2,
B^n_2)]^n=[as_{\phi_1}(B^n_2)]^{n-k} [as_{\phi_2}(B^n_2)]^k$,
\begin{eqnarray*}\left[\frac{as_k(\phi_1, K; \phi_2,
B^n_2)}{as_k(\phi_1, B^n_2; \phi_2, B^n_2)}\right]^n &\geq&
\left[\frac{as_{\phi_1}(K)}{as_{\phi_1}(B^n_2)}\right]^{n-k} \geq
\left(\frac{|K|}{|B^n_2|}\right)^{(n-k)(1-2r_1)},\end{eqnarray*}
where the last inequality follows inequality
(\ref{affine:isoperimetric:inequality:homo}) and $n-k\leq 0$.
Clearly if $K$ is a ball, the equality holds.

\vskip 1mm \noindent Theorem \ref{Stantalo:phi} implies that
$as_{\phi_1}(K)as_{\phi_1}(K^\circ)\leq [as_{\phi_1}(B^n_2)]^2$.
Combining with inequality (\ref{i:mixed:phi:2}) and $n-k\leq 0$,
we have
\begin{eqnarray*}
as_{k}(\phi_1, K; \phi_2, B^n_2)as_{k}(\phi_1, K^\circ; \phi_2,
B^n_2)&\geq&
[as_{\phi_1}(K)as_{\phi_1}(K^\circ)]^{\frac{n-k}{n}}[as_{\phi_2}(B^n_2)]^{\frac{2k}{n}}\\
&\geq& [as_{\phi_1}(B^n_2)]^{\frac{2(n-k)}{n}}[as_{\phi_2}(B^n_2)]^{\frac{2k}{n}}\\
&=&[as_{k}(\phi_1, B^n_2; \phi_2, B^n_2)]^2.
\end{eqnarray*} Clearly if $K$ is a ball, equality holds true.

\begin{proposition} \label{ith:mixed:general:k:smaller:0} Let $K\in C^2_+$ be a convex body with
centroid at the origin. For $k\leq 0$, and $\psi_1, \psi_2 \in
Conv(0, \infty)$, one has

\vskip 1mm \noindent (i) $as_k(\psi_1,\! K;\! \psi_2,\!
B^n_2)\!\geq\! as_k(\psi_1,\! B_K;\! \psi_2,\! B^n_2),$ with
equality if $K$ is a ball;

\vskip 1mm \noindent (ii) if in addition, $\psi_1$ is homogeneous
of degree $r_1\in (-\infty, 0]$, then
$$\left(\frac{as_k(\psi_1, K; \psi_2, B^n_2)}{as_k(\psi_1, B^n_2; \psi_2, B^n_2)}\right)^n\geq
\left(\frac{|K|}{|B^n_2|}\right)^{(n-k)(1-2r_1)},$$ with equality
if $K$ is a ball; Moreover, $$as_k(\psi_1, K; \psi_2,
B^n_2)as_k(\psi_1, K^\circ; \psi_2, B^n_2)\geq c^{n-k}
[as_k(\psi_1, B^n_2; \psi_2, B^n_2)]^2,$$ where $c$ is the
universal constant in the inverse Santal\'{o} inequality \cite{BM,
GK2, MilmanPajor2000, Nazarov2009, Pisier1989}; namely,
$|K||K^\circ|\geq c^n|B^n_2|^2$.
\end{proposition}

\vskip 1mm \noindent {\bf Proof.} (i) For the $L_{\psi}$-affine
surface area with $\psi \!\in\! Conv(0,\infty)$, Ludwig proved the
affine isoperimetric inequality \cite{Ludwig2009}:
$as_{\psi}(K)\!\geq\! as_{\psi}(B_K)$ with equality if and only if
$K$ is an ellipsoid. If $\psi\!\in\! Conv(0, \infty)$ is
homogeneous of degree $r\!\in\! (-\infty, 0]$, then
\begin{eqnarray}\label{affine:isoperimetric:inequality:homo:psi}\frac{as_{\psi}(K)}{as_{\psi}(B^n_2)}\geq
\left(\frac{|K|}{|B^n_2|}\right)^{1-2r},\end{eqnarray} with
equality if and only if $K$ is an ellipsoid.

\vskip 1mm  Similar to inequality (\ref{i:mixed:phi:2}), one can
prove that, for $k\leq 0$,
\begin{equation}\label{i:mixed:psi:2}\big[as_k(\psi_1, K; \psi_2,
L)\big]^n\geq
[as_{\psi_1}(K)]^{n-k}[as_{\psi_2}(L)]^{k}.\end{equation}
Combining with $[as_k(\psi_1, B_K; \psi_2,
B^n_2)]^n=[as_{\psi_1}(B_K)]^{n-k} [as_{\psi_2}(B^n_2)]^k$, one
has
\begin{eqnarray*}\left[\frac{as_k(\psi_1, K; \psi_2,
B^n_2)}{as_k(\psi_1, B_K; \psi_2, B^n_2)}\right]^n &\geq&
\left[\frac{as_{\psi_1}(K)}{as_{\psi_1}(B_K)}\right]^{n-k}\geq
1,\end{eqnarray*} where we have used Ludwig's isoperimetric
inequality in \cite{Ludwig2009} and $k\leq 0$. Equality holds
trivially if $K$ is a ball.

\vskip 1mm \noindent (ii)  Again by inequality
(\ref{i:mixed:psi:2}) and $[as_k(\psi_1, B^n_2; \psi_2,
B^n_2)]^n=[as_{\psi_1}(B^n_2)]^{n-k} [as_{\psi_2}(B^n_2)]^k$,
\begin{eqnarray*}\left[\frac{as_k(\psi_1, K; \psi_2,
B^n_2)}{as_k(\psi_1, B^n_2; \psi_2, B^n_2)}\right]^n &\geq&
\left[\frac{as_{\psi_1}(K)}{as_{\psi_1}(B^n_2)}\right]^{n-k}\geq
\left(\frac{|K|}{|B^n_2|}\right)^{(n-k)(1-2r_1)},\end{eqnarray*}
where the last inequality follows inequality
(\ref{affine:isoperimetric:inequality:homo:psi}) and $n-k\geq 0$.
Clearly if $K$ is a ball, the equality holds.

\vskip 1mm The following inequality was proved in
\cite{Ludwig2009}: if $\psi\in Conv(0, \infty)$, then
$$as_{\psi}(K)\geq n|K|\psi \left(\frac{|K^\circ|}{|K|}\right).$$
(Note that it holds true for all $K\in C^2_+$). The inverse
Santal\'{o} inequality says that $|K||K^\circ|\geq c^n |B^n_2|^2$,
where $c$ is a universal constant \cite{BM, GK2, MilmanPajor2000,
Nazarov2009, Pisier1989}. Therefore
$$as_{\psi}(K)as_{\psi}(K^\circ)\geq n^2|K||K^\circ|\psi (1)^2\geq
c^n\psi(1)^2 n^2|B^n_2|^2=c^n\big[as_{\psi}(B^n_2)\big]^2.$$
Combining with inequality (\ref{i:mixed:psi:2}) and $n-k\geq 0$,
\begin{eqnarray*}
as_{k}(\psi_1, K; \psi_2, B^n_2)as_{k}(\psi_1, K^\circ; \psi_2,
B^n_2)&\geq&
[as_{\psi_1}(K)as_{\psi_1}(K^\circ)]^{\frac{n-k}{n}}[as_{\psi_2}(B^n_2)]^{\frac{2k}{n}}\\
&\geq& c^{n-k}[ as_{\psi_1}(B^n_2)]^{\frac{2(n-k)}{n}}[as_{\psi_2}(B^n_2)]^{\frac{2k}{n}}\\
&=&c^{n-k}[as_{k}(\psi_1, B^n_2; \psi_2, B^n_2)]^2.
\end{eqnarray*}

\begin{proposition}
\label{ith:mixed:general:k:smaller:0:star} Let $K\in C^2_+$ be a
convex body with centroid at the origin. For $k\leq 0$, and
$\psi_1, \psi_2 \in Conv(0, \infty)$, one has

\vskip 1mm \noindent (i) $as_k^*(\psi_1,\! K;\! \psi_2,\!
B^n_2)\!\geq\! as_k^*(\psi_1,\! (B_{K^\circ})^\circ;\! \psi_2,\!
B^n_2),$ with equality if $K$ is a ball;

\vskip 1mm \noindent (ii) if in addition, $\psi_1$ is homogeneous
of degree $r_1\in (-\infty, 0]$, one has
$$\left(\frac{as_k^*(\psi_1, K; \psi_2, B^n_2)}{as_k^*(\psi_1, B^n_2; \psi_2,
B^n_2)}\right)^n\!\geq\!
c^{n(1-2r_1)(n-k)}\left(\frac{|K|}{|B^n_2|}\right)^{(n-k)(2r_1-1)},$$
Moreover, $as_k^*(\psi_1, K; \psi_2, B^n_2)as_k^*(\psi_1, K^\circ;
\psi_2, B^n_2)\geq c^{n-k} [as_k^*(\psi_1, B^n_2; \psi_2,
B^n_2)]^2,$ where $c$ is the same constant as in Proposition
\ref{ith:mixed:general:k:smaller:0}.
\end{proposition}

\vskip 1mm \noindent {\bf Proof.} (i) Recall that
$as_{\psi}^*(K)=as_{\psi}(K^\circ)$ (see Proposition
\ref{Duality:psi}). By inequality
(\ref{affine:isoperimetric:inequality:homo:psi}), and Ludwig's
isoperimetric inequality in \cite{Ludwig2009},
\begin{eqnarray*}\left[\frac{as^*_k(\psi_1,\! K;\! \psi_2,\!
B^n_2)}{as^*_k(\psi_1,\! (B_{K^\circ})^\circ;\! \psi_2,\!
B^n_2)}\right]^n \!\geq\!
\left[\frac{as^*_{\psi_1}(K)}{as^*_{\psi_1}((B_{K^\circ})^\circ)}\right]^{n-k}\!=\!
\left[\frac{as_{\psi_1}(K^\circ)}{as_{\psi_1}(B_{K^\circ})}\right]^{n-k}\!\geq\!
1.\end{eqnarray*} Clearly, the equality holds if $K$ is a ball.

\vskip 1mm \noindent (ii) Again, by inequality
(\ref{affine:isoperimetric:inequality:homo:psi}), the inverse
Santal\'{o} inequality and $r_1\in (-\infty, 0]$,
\begin{equation}\label{psi:santalo:1}\frac{as_{\psi}^*(K)}{as_{\psi}^*(B^n_2)}=\frac{as_{\psi}(K^\circ)}{as_{\psi}(B^n_2)}\geq
\left(\frac{|K^\circ|}{|B^n_2|}\right)^{1-2r_1} \geq
c^{n(1-2r_1)}\left(\frac{|K|}{|B^n_2|}\right)^{2r_1-1}.\end{equation}
Similar to inequality (\ref{i:mixed:phi:2}), one can prove that
for $k\leq 0$,
\begin{equation*}\big[as_k^*(\psi_1, K; \psi_2,
L)\big]^n\geq
[as^*_{\psi_1}(K)]^{n-k}[as^*_{\psi_2}(L)]^{k}.\end{equation*} As
$n-k\geq 0$ and inequality (\ref{psi:santalo:1}), one has
\begin{eqnarray*}\left[\frac{as^*_k(\psi_1, K; \psi_2,
B^n_2)}{as^*_k(\psi_1, B^n_2; \psi_2, B^n_2)}\right]^n &\geq&
\left[\frac{as^*_{\psi_1}(K)}{as^*_{\psi_1}(B^n_2)}\right]^{n-k}\geq
c^{n(1-2r_1)(n-k)}\left(\frac{|K|}{|B^n_2|}\right)^{(2r_1-1)(n-k)}.\end{eqnarray*}
The proof of $as_k(\psi_1, K; \psi_2, B^n_2)as_k(\psi_1, K^\circ;
\psi_2, B^n_2)\geq c^{n-k} [as_k(\psi_1, B^n_2; \psi_2, B^n_2)]^2$
is same as that of Proposition
\ref{ith:mixed:general:k:smaller:0}.

 \vskip 1mm
\noindent {\bf Acknowledgments.} This paper was completed by the
support of NSF-FRG DMS: 0652571 and 0652684, during the author's
postdoctoral fellowship at University of Missouri, Columbia. The
author is grateful to Dr.\ Werner and the reviewer for their
valuable comments.

\vskip 2mm \noindent Deping Ye \\
{\small Department of Mathematics and Statistics}\\
{\small Memorial University of Newfoundland}\\
{\small St. John's, Newfoundland, Canada A1C 5S7.} \\ {\tt
deping.ye@mun.ca}

\end{document}